\documentclass[pdflatex,sn-mathphys-num,iicol]{sn-jnl}

\usepackage{graphicx}%
\usepackage{multirow}%
\usepackage{amsmath,amsfonts}%
\usepackage{amsthm}%
\usepackage{amssymb}%
\usepackage{mathtools}   
\usepackage{fix-cm}
\usepackage{mathrsfs}%
\usepackage[title]{appendix}%
\usepackage{xcolor}%
\usepackage{textcomp}%
\usepackage{manyfoot}%
\usepackage{booktabs}%
\usepackage{algorithm}%
\usepackage{algorithmicx}%
\usepackage{algpseudocode}%
\usepackage{listings}%

\usepackage{color}
\usepackage{float}
\usepackage{url}
\usepackage{mdframed}
\usepackage{changes}
\usepackage{hyperref}
\usepackage{comment}
\usepackage{makecell}
\usepackage{anyfontsize}

\usepackage{subfigure} 
\usepackage{caption} 
\usepackage{multirow}
\usepackage{multicol}    
\usepackage{array}       
\usepackage{graphicx}
\graphicspath{{Figs/}}
\usepackage{empheq}
\allowdisplaybreaks
\theoremstyle{thmstyleone}%
\theoremstyle{thmstyletwo}%

\theoremstyle{thmstylethree}%

\newcommand{\bi}{\begin{itemize}}
\newcommand{\ei}{\end{itemize}}
\newcommand{\ba}{\begin{array}}
\newcommand{\ea}{\end{array}}

\begin{document}

\title[Article Title]{Robust Sparse Identification of Nonlinear Dynamics via Least Trimmed Squares}


\author*[1]{\fnm{Fabio} \sur{Amaral}}\email{fabio.amaral@unesp.br}

\author[2]{\fnm{Geovani N.} \sur{Grapiglia}}\email{geovani.grapiglia@uclouvain.be}
\equalcont{These authors contributed equally to this work.}

\author[1]{\fnm{Cassio M.} \sur{Oishi}}\email{cassio.oishi@unesp.br}
\equalcont{These authors contributed equally to this work.}

\affil*[1]{\orgdiv{Departamento de Matemática e Computação}, \orgname{Universidade Estadual Paulista `Júlio de Mesquita Filho'}, \orgaddress{\street{R. Roberto Simonsen, 313}, \city{Presidente Prudente}, \postcode{19060-900}, \state{São Paulo}, \country{Brazil}}}

\affil[2]{\orgdiv{ICTEAM Institute}, \orgname{Universit\'e catholique de Louvain}, \orgaddress{\street{Avenue Georges Lema\^itre 4}, \city{Ottignies-Louvain-la-Neuve}, \postcode{1348}, \state{Walloon Brabant}, \country{Belgium}}}


\abstract{In this work, we propose a robust Sparse Identification of Nonlinear Dynamics (SINDy) pipeline for handling datasets corrupted by noise and outliers. The method decouples outlier filtering from sparse regression by combining Iterative Least Trimmed Squares (ILTS) with Sequentially Thresholded Least Squares (STLS). Unlike standard approaches that treat all observations uniformly within a single regression stage, the proposed ILTS-SINDy framework first applies an ILTS procedure that iteratively minimizes the sum of the smallest squared residuals to identify the most reliable observations without prior knowledge of outliers, after which STLS is used to recover a parsimonious governing model. Extensive numerical experiments across six benchmarks show that ILTS-SINDy achieves superior model discovery (structure) accuracy, is capable of exact structural recovery, and attains up to two orders of magnitude lower coefficient errors under severe corruption. Furthermore, ILTS-SINDy delivers these gains with modest computational overhead, operating an order of magnitude faster than ensemble and weak form approaches.}




\maketitle
\section{Introduction}\label{sec:1}

Among the various classes of mathematical models, differential equations provide a fundamental framework for describing the temporal evolution of dynamical systems. However, deriving governing equations from first principles is often a challenging task that requires detailed knowledge of the phenomenon under study. This difficulty has motivated, in recent years, the development of data-driven approaches that aim to infer differential equations directly from observational data. In this context, the Sparse Identification of Nonlinear Dynamics (SINDy) framework, introduced by Brunton, Proctor, and Kutz~\cite{Brunton:2016}, has emerged as a particularly influential method for discovering parsimonious dynamical models from data. Given measurements of the system states $\{x_i(t_\ell)\}_{\ell=0}^m$, $i=1,\ldots,N$, SINDy assumes that the underlying dynamics can be represented as
\begin{equation}
\dot{x}_i(t) = \sum_{j=1}^{d} \alpha_{i,j}\,\theta_j(x(t)), \qquad i=1,\ldots,N,
\label{eq:sindy_problem}
\end{equation}
where $x(t) \in \mathbb{R}^N$, and $\{\theta_j\}_{j=1}^d$ is a user-defined library of candidate functions $\theta_j:\mathbb{R}^{N} \to \mathbb{R}$. Evaluating the system at the measurement points and approximating the time derivatives $\dot{x}_i(t_\ell)$ using finite-difference schemes, the coefficients $\alpha_{i,j}$ are estimated by solving $N$ linear regression problems via least squares, often combined with sparsity-promoting techniques to identify a low-dimensional and interpretable representation of the governing equations.

SINDy has been successfully applied across a broad range of scientific domains, including macroeconomic systems~\cite{LI}, oscillatory biological dynamics~\cite{Prokop}, and epidemiological modeling of COVID-19 transmission~\cite{Jiang2}. In combination with proper orthogonal decomposition, it has also been used for the identification of parametric models in viscoelastic fluid dynamics~\cite{Oishi}. Despite these successes, the standard SINDy framework is known to be highly sensitive to measurement noise and, in particular, to outliers. Since derivative estimation relies on finite differences of noisy observations, corrupted data can significantly distort derivative approximations, which in turn deteriorates the quality of the regression problem. This contamination may lead to severely biased models, thereby affecting both interpretability and predictive performance (see, e.g.,~\cite{Prokop,LI}). This limitation is especially critical in real-world applications, where data are typically noisy and may contain significant outliers \cite{Singh2025}.

To address these challenges, several extensions of SINDy have been proposed to improve robustness. A first class of methods explicitly targets outlier removal through data selection strategies. In this direction, trimmed SINDy approaches \cite{Champion:2020,kiser:2023} remove some data points that fit the model least well during regression, so that the model is learned from more reliable data. Another important development is Weak SINDy (W-SINDy)~\cite{Messenger:2021}, which reformulates the identification problem in a weak (integral) form. By applying integration by parts, derivatives are transferred onto smooth test functions, reducing the reliance on numerical differentiation and mitigating noise amplification. Ensemble SINDy (E-SINDy)~\cite{Fasel:2022} further improves robustness by training multiple models on different resampled versions of the data and combining their results, while also providing a natural way to quantify uncertainty.

Despite these advances, each approach has inherent limitations. W-SINDy improves robustness to measurement noise but introduces additional methodological complexity and requires careful selection of test functions and integration domains, which may limit its ease of use in practice. E-SINDy enhances stability and provides uncertainty estimates, but it can be computationally expensive and may still inherit biases present in the individual resampled datasets, particularly under severe data corruption. 

\begin{figure*}[htb!]
    \centering
    \includegraphics[width=\linewidth]{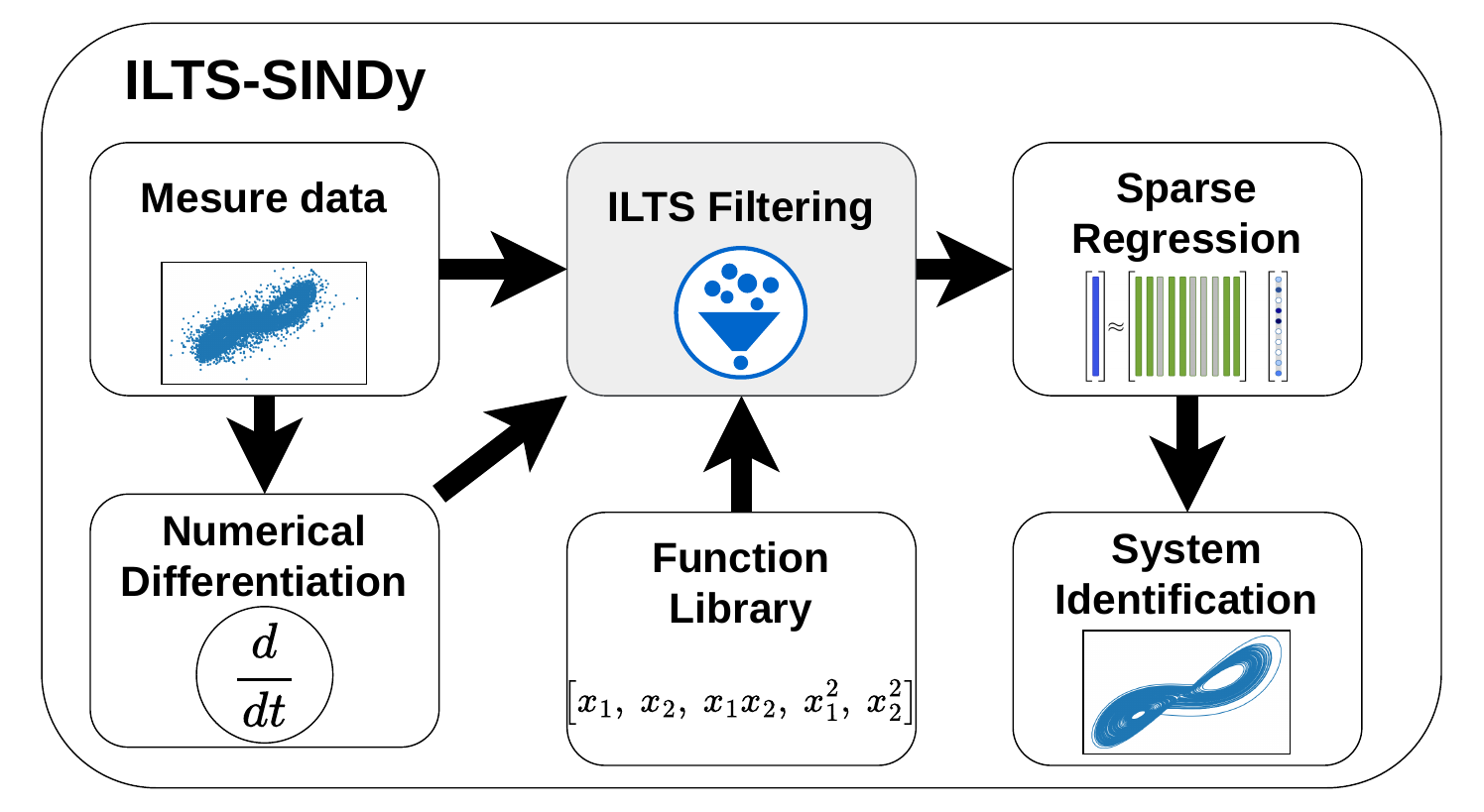}
    \caption{Overview of the proposed ILTS-SINDy framework for robust model discovery. Measured trajectories are differentiated numerically, the samples most affected by outliers are filtered out, and a sparse regression is performed on the remaining data to recover the governing equations.}
    \label{fig:method_overview}
\end{figure*}

Trimmed SINDy methods explicitly account for corrupted observations, but the identification of the reliable data subset and the estimation of the sparse governing equations are typically performed jointly. Consequently, data selection and sparse model identification remain coupled, making it difficult to assess their individual effects and potentially increasing the sensitivity of the identification procedure to the chosen optimization and sparsity parameters. These observations motivate a strategy in which corrupted observations are identified and removed in a dedicated preprocessing stage before performing sparse model identification.
In this context, the contributions of the present work are twofold:
\\[0.2cm]
\textbf{1. A robust filter-then-sparsify framework:} We revisit the trimmed SINDy framework and propose a modified identification pipeline in which data filtering and sparse model identification are decoupled into two sequential stages. Specifically, we first apply an iterative trimming procedure to detect and remove corrupted or highly inconsistent observations. For the data filtering stage, we employ Iterative Least Trimmed Squares (ILTS) \cite{rousseeuw:2006,shen2019iterative}, which alternates between estimating regression coefficients and retaining observations with the smallest residuals, thereby concentrating on the most consistent data points and mitigating the influence of outliers without requiring prior knowledge of their locations. For the sparsity-inducing stage, we use Sequentially Thresholded Least Squares (STLS), which iteratively thresholds small coefficients and recomputes the regression using only the remaining active terms, progressively refining the identified model until convergence to a potentially parsimonious representation of the governing equations. The resulting pipeline, referred to as ILTS-SINDy, is summarized in Figure~\ref{fig:method_overview}.
\vspace{0.3cm}

\noindent\textbf{2. A systematic numerical evaluation:} We evaluate the proposed ILTS-SINDy framework through extensive numerical experiments on synthetic benchmark datasets obtained from the the Susceptible–-Infected-–Recovered (SIR) model \cite{kermack1927contribution}, Lorenz \cite{lorenz1963deterministic},  Rössler \cite{rossler:1976},  ABC (Arnold-Beltrami-Childress) flow \cite{Gromeka:1881, Viknesh:2026}, Van der Pol oscillator \cite{vanderpol:1926} and the Lotka–Volterra \cite{lotka1925elements,volterra1926variazioni} dynamical systems. In particular, we investigate the robustness of the ILTS-SINDy filtering stage under varying levels of noise and data corruption. The results show a consistent improvement over standard SINDy with STLS, E-SINDy and W-SINDy, particularly in regimes of strong data contamination.

The remainder of this paper is organized as follows. Section 2 introduces the mathematical formulation of SINDy, some of its robust variants, and our proposed ILTS-SINDy pipeline. Section~3 presents numerical experiments on synthetic benchmark systems and provides a comparative evaluation against classical SINDy variants. Finally, Section~4 concludes the paper with a discussion of the main findings and directions for future work.
\\[0.2cm]
\noindent\textbf{Notation.} Given a natural number $m\in\mathbb{N}$, we denote $[m]:=\left\{0,1,\ldots,m\right\}$. The cardinality of a set $\mathcal{A}$ is denoted by $|\mathcal{A}|$.
\vspace{-0.2cm}
\section{Background and Methods}
\vspace{-0.2cm}
In this section, we start with a brief description about SINDy and some of its robust variants. Then we provide the details of our new approach, ILTS-SINDy.

\subsection{The SINDy Framework}

Using the notation introduced in the Introduction, let
\[
x(t_\ell)=\left[x_1(t_\ell),\ldots,x_N(t_\ell)\right]^T,
\qquad \ell=0,\ldots,m,
\]
denote the measured state vectors. For each measurement point $t_\ell$, the candidate functions
$\{\theta_j\}_{j=1}^d$ are evaluated to construct the library matrix
\[
A =
\begin{bmatrix}
\theta_1(x(t_0)) & \cdots & \theta_d(x(t_0))\\
\theta_1(x(t_1)) & \cdots & \theta_d(x(t_1))\\
\vdots           & \ddots & \vdots\\
\theta_1(x(t_m)) & \cdots & \theta_d(x(t_m))
\end{bmatrix}
\in\mathbb{R}^{(m+1)\times d}.
\]
The time derivatives $\dot{x}_{i}(t_{\ell})$ are approximated from the measured data using central finite-differences. Assuming uniformly sampled measurements with time step $\Delta t=t_{\ell+1}-t_\ell$, the entries of the observation vector $b_i\in\mathbb{R}^{m+1}$ approximating the derivatives $\dot{x}_{i}(t_{\ell})$, $1,\ldots,m$, are defined by
\[
b_i =
\begin{bmatrix}
b_{i,0}&
b_{i,1}&
\ldots&
b_{i,m}
\end{bmatrix}^{T},
\]
with
\[
b_{i\ell}=
\begin{cases}
\dfrac{x_i(t_{1})-x_i(t_0)}{\Delta t},
& \ell=0,\\[2ex]
\dfrac{x_{i}(t_{\ell+1})-x_{i}(t_{\ell-1})}{2\Delta t}, & \ell=1,\ldots,m-1,\\[2ex]
\dfrac{x_i(t_\ell)-x_i(t_{\ell-1})}{\Delta t},
& \ell=m.
\end{cases}
\]
Then, for each target state variable $x_i$, the SINDy problem for identifying the dynamics in (\ref{eq:sindy_problem}) is reduced to the linear least-squares problem
\begin{equation}
\min_{\alpha_i\in\mathbb{R}^d}
\frac12 \|A\alpha_i-b_i\|_2^2,
\label{eq:sindy_LS}
\end{equation}
under the assumption that the coefficient vector $\alpha_i=\left[\begin{array}{ccc}\alpha_{i,1}&\ldots &\alpha_{i,d}\end{array}\right]^{T}$ is sparse. Since only a small subset of candidate functions is expected to contribute to the true dynamics, sparsity-promoting methods are employed to eliminate insignificant coefficients.

A commonly used approach in SINDy is the STLS algorithm. Starting from the ordinary least-squares solution, STLS iteratively refines the coefficient estimate by alternating between hard thresholding and restricted least-squares regression. At each iteration, coefficients with magnitude below a prescribed threshold are set to zero, defining an active set of variables. The least-squares problem is then re-solved using only the remaining active variables, yielding updated coefficient values. These two steps are repeated until the active set no longer changes, resulting in a sparse approximate solution to the original linear regression problem. More precisely, a detailed description of the steps of STLS applied to (\ref{eq:sindy_LS}) is given in Algorithm \ref{alg:stls}.

\begin{algorithm}
\caption{STLS}
\label{alg:stls}
\begin{algorithmic}[1]
\Require Library matrix $A \in \mathbb{R}^{(m+1) \times d}$, target vector $b_i \in \mathbb{R}^{m+1}$, threshold parameter $\lambda > 0$
\Ensure $\tilde{\alpha}_{i}$

\State \textbf{Step 0 (Initialization).} Set initial active set $\mathcal{A}^{(0)} = \{1,\ldots,d\}$ and $k := 0$

\While{true}
    \State \textbf{Step 1 (Selection matrix construction).} Let $\mathcal{A}^{(k)} = \{j_1,\ldots,j_s\}$ with $s=|\mathcal{A}^{(k)}|$ and $j_1 < \cdots < j_s$
    \State Define selection matrix $U^{(k)} = [e_{j_1}, e_{j_2}, \ldots, e_{j_s}] \in \mathbb{R}^{d \times s}$
    
    \State \textbf{Step 2 (Restricted least-squares estimation).} Compute
    \State \quad $\beta_i^{(k)} = \arg\min_{\beta \in \mathbb{R}^s} \|A U^{(k)} \beta - b_i\|_2^2$
    \State \quad Lift to full space: $\alpha_i^{(k)} = U^{(k)} \beta_i^{(k)}$
    
    \State \textbf{Step 3 (Thresholding).} Update active set:
    \State \quad $\mathcal{A}^{(k+1)} = \left\{j \in \{1,\ldots,d\} : |\alpha_{i,j}^{(k)}| > \lambda\right\}$
    
    \State \textbf{Step 4 (Convergence check).} 
    \If{$\mathcal{A}^{(k+1)} = \mathcal{A}^{(k)}$}
        \State \Return $\tilde{\alpha}_{i}:=\alpha_i^{(k)}$
    \EndIf
    \State $k := k+1$
\EndWhile
\end{algorithmic}
\end{algorithm}

Although SINDy combined with STLS has demonstrated strong performance in recovering compact governing equations, its accuracy can be significantly affected by noise in both the data matrix $A$ and the estimated time derivatives used in the regression. In particular, measurement noise propagates through the numerical differentiation procedure, leading to errors in the target vector, while noise in the state measurements directly corrupts the construction of the regression matrix. These sources of error can degrade the quality of the identified model, motivating the development of more robust variants.

\subsection{SINDy with Trimming}

SINDy with trimming of outliers, introduced by \cite{Champion:2020}, extends standard sparse regression by explicitly accounting for corrupted measurements through a joint estimation of the model coefficients and an inlier selection mechanism. For identifying the Ordinary Differential Equation (ODE) describing $\dot{x}_{i}(t)$, this approach can be formulated as a linearly constrained composite optimization problem
\begin{equation}
    \begin{array}{l}
    \min\limits_{\substack{\alpha_i \in \mathbb{R}^d,\\ v \in \mathbb{R}^{m+1}}}
    \frac{1}{2} \sum_{\ell=0}^{m}
    v_\ell \left( \langle A_{\ell}, \alpha_i \rangle - b_{i\ell} \right)^2
    + \gamma \mathcal{R}(\alpha_i)\\\\
     \text{subject to}\\
    0 \leq v_\ell \leq 1,\quad \ell=0,\ldots,m,
    \quad \text{and}\quad
    \sum_{\ell=0}^{m} v_\ell = p,
    \end{array}
    \label{eq:sindy_trimmed}
\end{equation}
where \(A_{\ell} = [\theta_1(x(t_\ell)), \ldots, \theta_d(x(t_\ell))]^T\) denotes the \((\ell+1)\)-th row of the regression matrix, \(b_{i\ell}\) is the corresponding time-derivative of the observation for the \(i\)-th state component, \(\gamma>0\) is the regularization parameter and \(\mathcal{R}\) is a regularization function, such as \(\ell_0\) or \(\ell_1\)--norm. The auxiliary variables \(v_\ell\) act as soft inlier indicators: values close to one indicate that the corresponding data point is treated as an inlier, while values close to zero effectively downweight or discard its contribution in the regression. The constraint \(\sum_{\ell=0}^m v_\ell = p\) enforces that only \(p < m\) samples are selected as inliers, thereby inducing an implicit trimming of the dataset.

The resulting objective couples sparse regression with sample selection: the first term enforces a good fit on the selected subset of data, while the regularization term promotes sparsity in the coefficients, leading to a parsimonious representation of the dynamics. In practice, problem (\ref{eq:sindy_trimmed}) is typically solved via an alternating minimization strategy based on the Sparse Relaxed Regularized Regression (SR3) algorithm \cite{Champion}. This type of trimmed formulation introduces a fundamental trade-off between identifying corrupted observations and estimating sparse governing equations. Since both tasks are performed simultaneously, highly corrupted samples may still be classified as inliers, which can in turn degrade the quality of the recovered model.

\subsection{Ensemble SINDy}

E-SINDy~\cite{Fasel:2022} improves robustness of sparse regression by aggregating multiple independent model estimates obtained from resampled versions of the data. Instead of fitting a single model, E-SINDy constructs an ensemble of sparse solutions, each computed from a different subsample of the available measurements. For each realization \(r = 1,\ldots,K\), a sub-sampled regression problem is formed and a sparse coefficient vector $\alpha_i^{(r)} =
\begin{bmatrix}
\alpha_{i1}^{(r)} & \cdots & \alpha_{id}^{(r)}
\end{bmatrix}^T \in \mathbb{R}^d
$ is obtained using a procedure such as STLS {(Algorithm \ref{alg:stls})}. This yields a collection of independent sparse models \(\{\alpha_i^{(r)}\}_{t=1}^K\), each reflecting the variability induced by data resampling. To identify consistently selected terms, the inclusion probability of the \(j\)-th coefficient is defined as
\[
p_j = \frac{1}{K}\sum_{t=1}^K \mathbf{1}\left(\alpha_{ij}^{(r)} \neq 0\right),
\]
and an active set is selected as $\mathcal{I} = \{ j : p_j \ge \tau \}$, for a prescribed threshold \(\tau \in [0,1]\). The final coefficient vector is then obtained through a robust component-wise aggregation of the ensemble. In the so-called \textit{bragging} formulation, coefficients are computed as the median across realizations and subsequently thresholded according to the selected active set:
\[
\alpha_{ij}^{\mathrm{final}} = \begin{cases}
\mathrm{median}\big(\alpha_{ij}^{(1)}, \ldots, \alpha_{ij}^{(K)}\big) ,\quad \text{for } j \in \mathcal{I};\\
0, \quad \text{otherwise}.
\end{cases}
\]

This ensemble strategy yields a robust sparse estimator by replacing a single potentially unstable regression with a consensus-based aggregation over multiple realizations, thereby reducing sensitivity to noise and sampling variability. Nevertheless, it can be computationally expensive and may still inherit biases present in individual samples, particularly under severe data corruption.

\subsection{Weak SINDy}

Numerical differentiation of noisy measurements, particularly by finite differences, amplifies noise and can severely degrade the regression step of SINDy. To circumvent this issue, Messenger and Bortz~\cite{Messenger:2021} proposed W-SINDy, a weak formulation of SINDy. The central idea is to multiply the governing equations by smooth, compactly supported test functions and integrate by parts, thereby transferring the time derivative from the noisy state measurements onto the test functions, whose derivatives are known analytically. In addition, the resulting integrals act as local averaging operators, which further attenuates measurement noise.

Consider an autonomous dynamical system whose $i$-th state component obeys
\begin{equation}
    \dot{x}_i(t) = \sum_{j=1}^{d} \alpha_{i,j}\, \theta_j({x}(t)), \qquad i = 1, \dots, N,
    \label{eq:wsindy_ode}
\end{equation}
where $\{\theta_j\}_{j=1}^{d}$ are the candidate library functions and $\alpha_{i,j}$ are the unknown coefficients. Let $\phi$ be a continuously differentiable test function supported on an interval $[T_0, T_f]$, such that $\phi(T_0) = \phi(T_f) = 0$. Multiplying Eq.~\eqref{eq:wsindy_ode} by $\phi$ and integrating by parts over $[T_0, T_f]$ yields
\begin{equation}
    -\int_{T_0}^{T_f} \phi'(t) x_i(t)\, dt = \sum_{j=1}^{d} \alpha_{i,j} \int_{T_0}^{T_f} \phi(t) \theta_j({x}(t))\, dt,
    \label{eq:wsindy_continuous}
\end{equation}
where the boundary term $\left[\phi(t)\, x_i(t)\right]_{T_0}^{T_f}$ vanishes because $\phi$ is zero at both endpoints. Note that Eq.~\eqref{eq:wsindy_continuous} no longer involves derivatives of the state.

In practice, the state is sampled on a uniform grid $\{t_k\}_{k=0}^{m}$ covering the observation window, and each test function is supported within this window. Approximating the integrals in Eq.~\eqref{eq:wsindy_continuous} with a quadrature rule (e.g., the trapezoidal rule) with weights $\{w_k\}_{k=0}^{m}$ gives
\begin{equation}
\begin{aligned}
    &-\sum_{k=0}^{m} w_k\, \phi'(t_k)\, x_i(t_k) \approx\\
    &\sum_{j=1}^{d} \alpha_{i,j} \sum_{k=0}^{m} w_k\, \phi(t_k)\, \theta_j({x}(t_k)).
\end{aligned}
    \label{eq:wsindy_discrete}
\end{equation}

Evaluating Eq.~\eqref{eq:wsindy_discrete} for a family of $K$ test functions $\{\phi_\ell\}_{\ell=1}^{K}$, each supported on a different subinterval, leads to the linear system ${G}\alpha_i \approx b^W_i$, which can be written compactly as
\begin{equation}
    {G} = {\Phi}\, {W}\, A
    \quad \text{and} \quad
    b^W_i = -{\Phi}'\, {W}\, {x}_i,
    \label{eq:wsindy_matrix_form}
\end{equation}
where
\begin{itemize}
    \item ${\Phi}, {\Phi}' \in \mathbb{R}^{K \times (m+1)}$ are the test function matrices with entries $[{\Phi}]_{\ell,k} = \phi_\ell(t_k)$ and $[{\Phi}']_{\ell,k} = \phi_\ell'(t_k)$;
    \item ${W} = \operatorname{diag}(w_0, \dots, w_m) \in \mathbb{R}^{(m+1) \times (m+1)}$ is the diagonal matrix of quadrature weights;
    \item $A \in \mathbb{R}^{(m+1) \times d}$ is the standard SINDy library evaluated at the sampling instants, with $[A]_{k,j} = \theta_j({x}(t_k))$;
    \item ${x}_i = [x_i(t_0), \dots, x_i(t_m)]^T \in \mathbb{R}^{m+1}$ is the sampled trajectory of the $i$-th state component.
\end{itemize}
Although Eq.~\eqref{eq:wsindy_matrix_form} plays the same role as the regression problem $A\alpha_i \approx b_i$ of standard SINDy, the two are constructed differently: the rows are indexed by the $K$ test functions rather than by the $m+1$ samples, so that $b_i^{\mathrm{W}} \in \mathbb{R}^{K}$, and each row summarizes the trajectory over a subinterval rather than at a single instant.

Since ${G}$ does not depend on $i$, it is assembled only once and shared by all state components. The coefficient vector $\alpha_i \in \mathbb{R}^{d}$ is then obtained by sparse regression,
\begin{equation}
    \alpha_i = \operatorname*{arg\,min}_{\alpha \in \mathbb{R}^{d}} \left\| {G}\alpha- b^W_i \right\|_2^2 + \lambda \mathcal{R}(\alpha),
    \label{eq:wsindy_regression}
\end{equation}
where $\mathcal{R}(\cdot)$ is a sparsity-promoting regularizer (or a thresholding criterion, as in STLS) and $\lambda > 0$ controls the level of sparsity.

A common choice of test functions is the family of piecewise polynomials
\begin{equation}
    \phi_\ell(t) =
    \begin{cases}
        C\, (t - T_\ell^{-})^{r}\, (T_\ell^{+} - t)^{s}, & t \in (T_\ell^{-}, T_\ell^{+}), \\
        0, & \text{otherwise},
    \end{cases}
    \label{eq:poly_test}
\end{equation}
where $C$ is a normalization constant (e.g., chosen such that $\max_t \phi_\ell(t) = 1$). The test functions are thus characterized by four hyperparameters. The exponents $r, s \geq 1$ control the smoothness and shape of each test function; in particular, $\phi_\ell$ is $(\min(r,s) - 1)$ times continuously differentiable, and $r = s$ yields a symmetric bell-shaped curve centered at $(T_\ell^{-} + T_\ell^{+})/2$. The support width $T_\ell^{+} - T_\ell^{-}$ sets the time scale over which the dynamics are averaged: wider supports provide stronger noise attenuation but may smooth out fast dynamics. Finally, the number of test functions $K$ determines the number of rows of the linear system; their centers are distributed over the observation window, either on a uniform grid or by random sampling, so that larger values of $K$ increase the overlap between adjacent supports. The values adopted in this work are reported in Section~\ref{sec:hyperparameters}.

\subsection{New Approach: ILTS-SINDy}
\vspace{-0.1cm}
Our starting point is the SINDy variant with trimming described in subsection 2.2. Restricting in (\ref{eq:sindy_trimmed}) the selection variables to binary values \(v_\ell \in \{0,1\}\) enforces a hard selection of exactly \(p\) observations through the constraint \(\sum_{\ell=0}^m v_\ell = p\). In this case, omitting the regularization term, the problem becomes equivalent to
\vspace{-0.2cm}
\begin{equation}
\min_{\alpha_{i}\subset\mathbb{R}^{d},\,\Omega\subset [m]}\,\frac{1}{2}\sum_{\ell\in\Omega} \left( \langle A_{\ell}, \alpha_i \rangle - b_{i,\ell} \right)^2\,\text{s.t}\, |\Omega|=p,  
\label{eq:LTS}
\end{equation}
which is the so called Least Trimmed Squares (LTS) problem. Defining the squared residuals $r_\ell(\alpha_i)^2 = \left( \langle A_{\ell}, \alpha_i \rangle - b_{i,\ell} \right)^2$, and ordering them as
\[
r_{\ell_{0}}(\alpha_i)^{2} \leq r_{\ell_{1}}(\alpha_i)^{2} \leq \cdots \leq r_{\ell_{m}}(\alpha_i)^{2},
\]
problem (\ref{eq:LTS}) can be written as
\begin{equation}
\min_{\alpha_i \in \mathbb{R}^d}
\frac{1}{2} \sum_{j=0}^{p-1} r_{\ell_{j}}(\alpha_i)^{2},
\label{eq:LOVO}
\end{equation}
which corresponds to selecting the \(p\) smallest residuals\footnote{In the Optimization community, problems of the form (\ref{eq:LOVO}) are also known as \textit{Low-Order Value Optimization} (LOVO) problems \cite{Andreani2}.}. To approximately solve this non-convex and non-smooth problem, one can use the ILTS method described in \cite{shen2019iterative}, and originally proposed in \cite{rousseeuw:2006}. Starting from an initial estimate of the coefficient vector, ILTS iteratively refines the solution by alternating between selecting the $p$ observations with the smallest residuals and solving a least-squares problem restricted to this subset. At each iteration, a new coefficient estimate is computed by re-fitting the model using the selected observations, and the procedure is repeated until the active set stabilizes, yielding an approximate solution of the original trimmed regression problem. In addition, ILTS returns the cleaned data pair $(\hat{A}, \hat{b}_i)$, where $\hat{A}\in\mathbb{R}^{p\times d}$ and $\hat{b}_i\in\mathbb{R}^{p}$ contain the rows and entries associated with the final active set, interpreted as the estimated uncorrupted samples. More precisely, a detailed description of ILTS applied to the LTS formulation is given in Algorithm \ref{alg:ilts}.

\begin{algorithm}
\caption{ILTS}
\label{alg:ilts}
\begin{algorithmic}[1]
\Require Matrix $A \in \mathbb{R}^{(m+1) \times d}$, target vector $b_i \in \mathbb{R}^{m+1}$, parameter $p < m$
\Ensure $\hat{\alpha}_{i}$, $\hat{A}$, $\hat{b}_{i}$

\State \textbf{Step 0 (Initialization).} Choose $\alpha_i^{(0)}\in\mathbb{R}^{d}$, set $\mathcal{S}^{(-1)}=\emptyset$, and $k := 0$

\While{true}
    \State \textbf{Step 1 (Residual computation and inlier selection).}
    \State Compute squared residuals: $r_\ell(\alpha_{i}^{(k)})^{2} = \left(\langle A_{\ell}, \alpha_i^{(k)} \rangle - b_{i,\ell}\right)^2$ for $\ell=0,\ldots,m$
    
    \State Order residuals as: $r_{\ell_{0}}(\alpha_{i}^{(k)})^{2} \leq \cdots \leq r_{\ell_{m}}(\alpha_{i}^{(k)})^{2}$
    
    \State Define $\mathcal{S}^{(k)}=\{\ell_0,\ldots,\ell_{p-1}\}$ (indices of $p$ smallest residuals)
    \State Set $r^{(k)} = \sum_{j = 0}^{p-1} r_{\ell_{j}}(\alpha_{i}^{(k)})^{2}$
    \State Define $\hat{A}^{(k)} = [A_{\ell_0}; \cdots; A_{\ell_{p-1}}]^T \in\mathbb{R}^{p\times d}$
    \State Define $\hat{b}_i^{(k)} = [b_{i,\ell_0}; \cdots; b_{i,\ell_{p-1}}]^T \in\mathbb{R}^{p}$
    
    \State \textbf{Step 2 (Stopping criterion).}
    \If{$k\geq 1$ and $r^{(k)}=r^{(k-1)}$}
        \State \Return $\hat{\alpha}_{i}:=\alpha_i^{(k)}$, $\hat{A}=\hat{A}^{(k)}$, $\hat{b}_{i}=\hat{b}_{i}^{(k)}$
    \EndIf
    
    \State \textbf{Step 3 (Least-squares update).}
    \State $\alpha_i^{(k+1)} = \arg\min_{\alpha_i \in \mathbb{R}^d} \frac{1}{2} \sum_{j=0}^{p-1} r_{\ell_{j}}(\alpha_i)^{2}$
    
    \State \textbf{Step 4 (Iteration update).}
    \State $k:=k+1$
\EndWhile
\end{algorithmic}
\end{algorithm}

To better understand the theoretical properties of Algorithm 2, it is useful to rewrite (\ref{eq:LOVO}) in an equivalent form. In particular, note that the number of subsets $\Omega\subset[m]$ with $|\Omega|=p$ is
\[
s=\binom{m}{p}.
\]
Enumerating these subsets as $\Omega_1,\ldots,\Omega_s$, problem (\ref{eq:LOVO}) can be written as
\[
\min_{\alpha_i\in\mathbb{R}^d}\varphi(\alpha_i),
\qquad
\varphi(\alpha_i):=\min_{j=1,\ldots,s}\varphi_j(\alpha_i),
\]
where
\[
\varphi_j(\alpha_i)
:=
\frac12\sum_{\ell\in\Omega_j}
\left(\langle A_\ell,\alpha_i\rangle-b_{i,\ell}\right)^2.
\]
We say that $\varphi_j(\,\cdot\,)$ is active at $\alpha_i$ if
\[
\varphi_j(\alpha_i)=\varphi(\alpha_i).
\]

In this notation, the steps of Algorithm 2 can be interpreted as follows. At the beginning of its $k$-th iteration, ILTS identifies an active function $\varphi_{j_k}(\,\cdot\,)$ at $\alpha_i^{(k)}$. It then computes a global minimizer of this function,
\[
\alpha_i^{(k+1)}
=
\arg\min_{\alpha_i\in\mathbb{R}^d}\varphi_{j_k}(\alpha_i).
\]
Consequently,
\[
\varphi(\alpha_i^{(k+1)})
\leq
\varphi_{j_k}(\alpha_i^{(k+1)})
\leq
\varphi_{j_k}(\alpha_i^{(k)})
=
\varphi(\alpha_i^{(k)}),
\]
so the objective function is non-increasing along the ILTS iterations. Moreover, as long as the objective value strictly decreases at each iteration, the same subset cannot be selected twice: selecting the same subset again would give the same least-squares minimizer and therefore could not produce a further decrease. Since there are $s=\binom{m}{p}$ possible subsets, ILTS can therefore perform at most $s$ strict-decrease iterations before  terminating. Thus, within a finite number of iterations, ILTS  returns a point that globally minimizes one of the functions defining the objective and for which that function $\varphi(\,\cdot\,)$ is active\footnote{This condition does not, in general, imply that the resulting point is a global minimizer of the original problem.}.

\begin{figure*}[htb!]
    \centering
    \includegraphics[width=\linewidth]{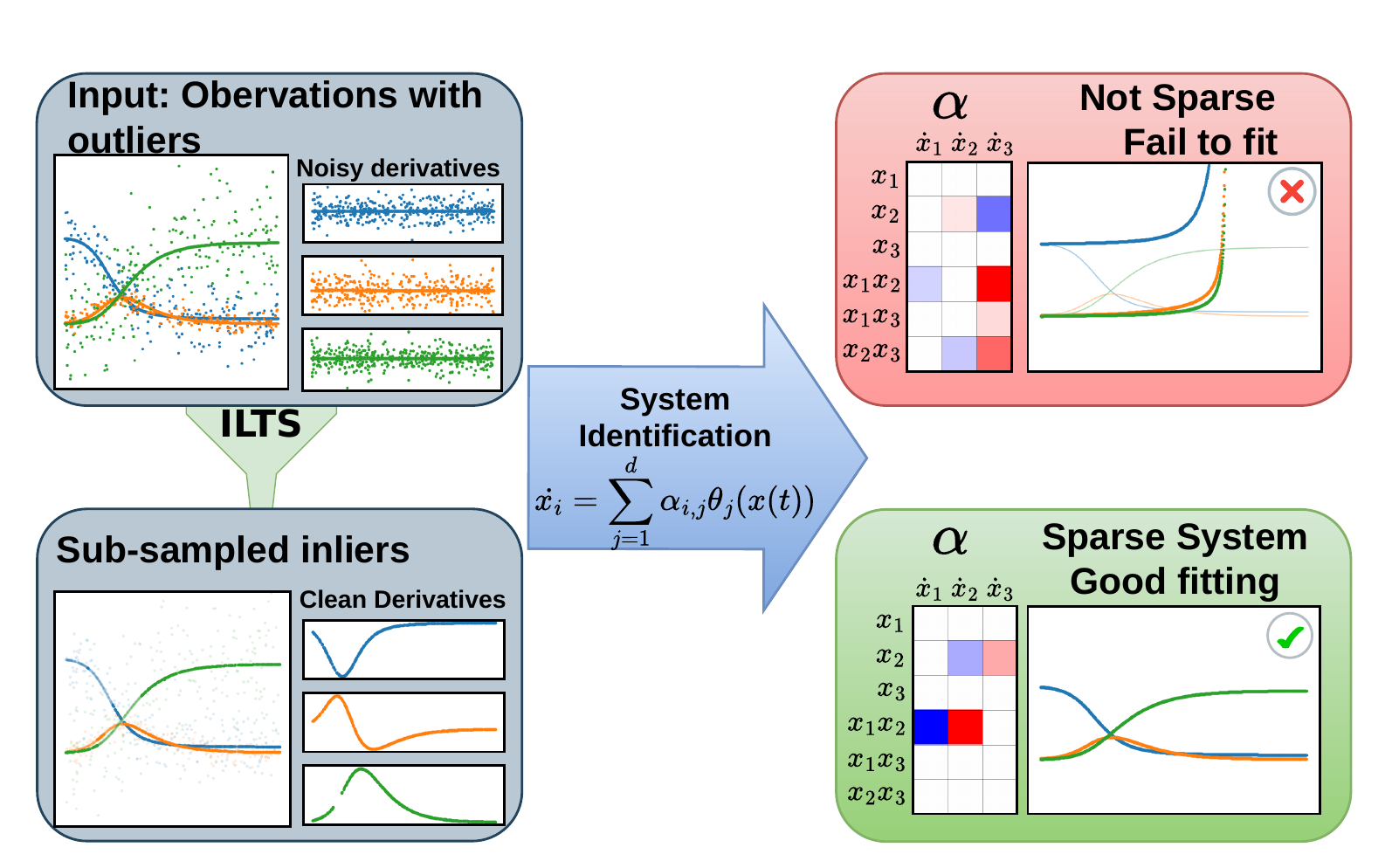}
    \caption{Standard SINDy approaches (top) are highly sensitive to contaminated data, often resulting in spurious model terms and divergent simulations. By contrast, ILTS-SINDy (bottom) automatically filters the data matrix to a subset of trusted inliers before sparse regression, allowing successful system identification even in the presence of outliers.}
    \label{fig:pipeline}
\end{figure*}

To achieve robust identification in the presence of severe data contamination, we introduce ILTS-SINDy, a simple pipeline in which corrupted samples are first filtered out using ILTS (Algorithm 2), after which STLS (Algorithm 1) is applied to the resulting cleaned data pair $(\hat{A}, \hat{b}_i)$ to recover a sparse coefficient vector $\alpha_i^{*}$ associated with the target state $x_i$. Repeating this procedure for all the states $x_{i}$, $i=1,\ldots,N$, we obtain the identified system of ODEs:

As mentioned above, standard SINDy approaches are highly sensitive to contaminated data since errors in both the regression matrix and the estimated time derivatives can propagate through the sparse regression step and often lead to spurious model terms. Trimmed SINDy formulations in \cite{Champion:2020,kiser:2023} mitigate this issue by adopting a \textit{filter-and-sparsify} strategy where outlier handling and sparse coefficient estimation are performed simultaneously within a single non-convex optimization problem. In contrast, ILTS-SINDy follows a \textit{filter-then-sparsify} strategy in which corrupted samples are first removed using ILTS and sparsification is subsequently carried out on the resulting cleaned dataset via STLS. This preliminary filtering step reduces the influence of outliers and produces a cleaner regression problem for sparse identification. Figure \ref{fig:pipeline} illustrates the proposed ILTS-SINDy approach in contrast to standard SINDy.

More precisely, ILTS-SINDy applied to the identification of the dynamics for the $i$-th state $x_{i}$ is described in Algorithm \ref{alg:ilts-sindy}.
\vspace{0.2cm}
\begin{algorithm}
\caption{ILTS-SINDy}
\label{alg:ilts-sindy}
\begin{algorithmic}[1]
\Require Matrix $A \in \mathbb{R}^{(m+1) \times d}$, target vector $b_i \in \mathbb{R}^{m+1}$, parameters $p < m$ and $\lambda>0$
\Ensure Sparse coefficient vector $\alpha_i^{*}$

\State \textbf{Step 0 (Initialization).}
\State Compute initial estimate: 
\State \quad $\alpha_i^{(0)} = \arg\min_{\alpha_{i}\in\mathbb{R}^{d}}\frac{1}{2}\|A\alpha_{i}-b_{i}\|_{2}^{2}$
\State Set $k := 0$

\State \textbf{Step 1 (Data trimming via ILTS).}
\State Apply Algorithm 2 to obtain cleaned data:
\State \quad $(\hat{\alpha}_{i},\hat{A},\hat{b}_{i})=\texttt{ILTS}(A,b_{i},p)$

\State \textbf{Step 2 (Sparse regression via STLS).}
\State Apply Algorithm 1 to cleaned data:
\State \quad $\alpha_i^{*} = \texttt{STLS}(\hat{A},\hat{b}_i,\lambda)$

\State \textbf{Step 3 (Output).}
\State \Return $\alpha_i^{*}$
\end{algorithmic}
\end{algorithm}

\vspace{0.3cm}
The convergence of the estimated coefficients to the true coefficients as the amount of data increases (i.e., consistency) requires correct identification of the active set, together with the usual assumptions for least-squares estimation. In this regard, STLS has limited theoretical guarantees, since the reliability of active-set selection depends on the identifiability and conditioning of the candidate library, the noise and derivative-estimation errors, the threshold parameter, and the magnitude of the nonzero coefficients. We therefore do not claim unconditional support recovery or consistency of the complete STLS procedure, but assess its robustness empirically, as described in the next section.

It worth noting that Algorithm~\ref{alg:ilts-sindy} take as input the hyperparamter $p$ and threshold $\lambda$. Here, $p$ is assumed to represent the number of inliers in the dataset, while $\lambda$ is the threshold used in STLS to enforce sparsity. In practice, applying Algorithm~\ref{alg:ilts-sindy} requiries estimating $p$ and tuning $\lambda$.

\section{Numerical Results}\label{sec:num_results}

This section presents a numerical comparison of ILTS-SINDy, standard SINDy,
W-SINDy, and E-SINDy. The comparison is based on synthetic
trajectories generated from six benchmark dynamical systems: the SIR model
\cite{kermack1927contribution}, the Lorenz system \cite{lorenz1963deterministic},
the R\"ossler system \cite{rossler:1976}, the ABC flow
\cite{Gromeka:1881, Viknesh:2026}, the Van der Pol oscillator
\cite{vanderpol:1926}, and the Lotka--Volterra system
\cite{lotka1925elements,volterra1926variazioni}. We first describe the models,
the procedure used to generate the synthetic noisy datasets, and the performance
metrics used in the numerical study.

Finally, we apply ILTS-SINDy to a real-world economic dataset and
compare it with the model discovered in \cite{li2025sparse}. The comparison is
carried out against the E-SINDy model reported in that reference, using both the
raw and the pre-processed data.

\subsection{Models, Data Generation and Performance Metrics}

The models, parameters, and initial conditions used for simulation are detailed in Table \ref{tab:models_params}. The SIR model describes the spread of an infectious disease within a population that, at time $t$, is partitioned into three compartments: susceptible individuals $x_1(t)$, infected individuals $x_2(t)$, and recovered individuals $x_3(t)$. The Lorenz model describes the evolution of a simplified atmospheric convection system, where $x_1(t)$ represents the intensity of convective motion, $x_2(t)$ corresponds to the horizontal temperature variation, and $x_3(t)$ denotes the vertical temperature variation. The R\"ossler system is a minimal chaotic model consisting of a linear spiral motion in the $(x_1, x_2)$ plane combined with a single quadratic nonlinearity in the equation for $x_3(t)$, which intermittently reinjects trajectories toward the plane and gives rise to a chaotic attractor of simpler geometry than that of the Lorenz system. The R\"ossler description gives no physical interpretation of the individual states. It was constructed as a mathematical prototype rather than derived from a physical system. The ABC flow is a steady, three-dimensional solution of the Euler equations in which the velocity and vorticity fields are everywhere parallel, and $x_1(t),\ x_2(t)$ and $x_3(t)$ denote the position of a fluid particle advected by the flow. Despite the steady and spatially periodic velocity field, particle trajectories are chaotic for generic parameter values. The Van der Pol oscillator is a non-conservative oscillating system with nonlinear damping, where $x_1(t)$ denotes the oscillator displacement and $x_2(t)$ its velocity. The damping is negative for small amplitudes and positive for large ones, so that trajectories converge to a self-sustained limit cycle. Finally, the Lotka--Volterra model describes the dynamics of a closed ecological system consisting of two interacting species, where $x_1(t)$ and $x_2(t)$ denote the prey and predator populations, respectively.

\small
\begin{table*}[htb]
    \caption{Benchmark models, parameter values, initial conditions, and simulation settings used in the numerical experiments. The time step is $\Delta t = t_{k+1} - t_{k}$, and $\lambda$ is the sparsity threshold of the STLS optimizer, applied identically to all methods.} \label{tab:models_params} 

    \begin{tabular}{|c|l|l|c|}
            \hline Model & \multicolumn{1}{c|}{Equations} & \multicolumn{1}{c|}{Parameters} & \multicolumn{1}{c|}{\begin{tabular}[c]{@{}c@{}}Initial\\Conditions\end{tabular}} \\ \hline SIR & 
            $\begin{aligned}
            \begin{cases}
                \dot{x}_1(t) &= -\beta x_1(t) x_2(t) \\
                \dot{x}_2(t) &= \beta x_1(t) x_2(t) - \gamma x_2(t) \\
                \dot{x}_3(t) &= \gamma x_2(t)
            \end{cases}
            \end{aligned}$& 
        \begin{tabular}[c]{@{}l@{}}
            $\beta = 0.3,\ \gamma = 0.1$\\
            $t \in [0,100],\ \Delta t = 0.1$\\
            threshold $\lambda:\ 0.01$
        \end{tabular} &
        $\left[\begin{array}{c} 0.99\\0.01\\0\end{array}\right]$ \\ \hline
        Lorenz & 
        $\begin{aligned}
            \begin{cases}
            \dot{x}_1(t) &= \sigma(x_2(t) - x_1(t)) \\
            \dot{x}_2(t) &= x_1(t)(\rho - x_3(t)) - x_2(t) \\
            \dot{x}_3(t) &= x_1(t)x_2(t) - \beta x_3(t)
            \end{cases}
        \end{aligned}$ 
        &
        \begin{tabular}[c]{@{}l@{}}
        $\sigma = 10, \beta = \frac{8}{3},\ \rho = 28$\\ $t \in [0,20],\ \Delta t = 0.02$\\
        threshold $\lambda:\ 0.1$
        \end{tabular} 
        & $\left[\begin{array}{c}-8\\7\\27\end{array}\right]$ \\ \hline
        R\"ossler & 
        $\begin{aligned}
            \begin{cases}
            \dot{x}_1(t) &= -x_2(t) - x_3(t) \\
            \dot{x}_2(t) &= x_1(t) + a x_2(t)\\
            \dot{x}_3(t) &= b + x_3(t) (x_1(t) - c) \\
            \end{cases}
        \end{aligned}$  & 
        \begin{tabular}[c]{@{}l@{}}
            $a = 0.2,\ b = 0.2,\ c = 5.7$\\
            $t \in [0,100],\ \Delta t = 0.1$\\
        threshold $\lambda:\ 0.1$
        \end{tabular}
        & $\left[\begin{array}{c} 0.1\\0.1\\ 0.1\end{array}\right]$ \\ 
        \hline 
        ABC Flow & 
        $\begin{aligned}
            \begin{cases}
            \dot{x}_1(t) &= A \sin(w x_3(t)) + C \cos(w x_2(t)) \\
            \dot{x}_2(t) &= B \sin(w x_1(t)) + A \cos(w x_3(t)) \\
            \dot{x}_3(t) &= C \sin(w x_2(t)) + B \cos(w x_1(t)) 
            \end{cases}
        \end{aligned}$  & 
        \begin{tabular}[c]{@{}l@{}}
            $A= 2,\ B = 3$,\\ $C = 1,\ w = 1$\\
            $t \in [0,20],\ \Delta t = 0.02$\\
            threshold $\lambda:\ 0.1$
        \end{tabular}
        & $\left[\begin{array}{c} 0.5\\0.2\\1\end{array}\right]$ \\ 
        \hline 
        Lotka-Volterra & 
        $\begin{aligned}
            \begin{cases}
            \dot{x}_1(t) &= \alpha x_1(t) - \beta x_1(t) x_2(t) \\
            \dot{x}_2(t) &= \beta x_1(t) x_2(t) - 2\alpha x_2(t)
            \end{cases}
        \end{aligned}$  & 
        \begin{tabular}[c]{@{}l@{}}
            $\alpha= 1,\ \beta = 0.1$\\
            $t \in [0,30],\ \Delta t = 0.03$\\
            threshold $\lambda:\ 0.05$
        \end{tabular}
        & $\left[\begin{array}{c} 1\\2\end{array}\right]$ \\ 
        \hline 
        Van Der Pol & 
        $\begin{aligned}
            \begin{cases}
            \dot{x}_1(t) &= x_2(t) \\
            \dot{x}_2(t) &= \mu x_2(t) (1 - x_1(t)^2) - x_1(t)
            \end{cases}
        \end{aligned}$  & 
        \begin{tabular}[c]{@{}l@{}}
            $\mu= 2$\\
            $t \in [0,100],\ \Delta t = 0.1$\\
            threshold $\lambda:\ 0.1$
        \end{tabular}
        & $\left[\begin{array}{c} 2\\0\end{array}\right]$ \\ 
        \hline 
    \end{tabular} 

\end{table*}

\normalsize

Each model, with the parameters specified in Table~\ref{tab:models_params}, was numerically integrated over the prescribed time interval using the adaptive explicit Runge--Kutta method of order 5(4) available in SciPy (\texttt{solve\_ivp} with method \texttt{RK45} and default tolerances), with the solution evaluated on a uniform grid with the time step $\Delta t$ given in Table~\ref{tab:models_params}. This yielded a clean trajectory, $\{\boldsymbol{x}(t_{\ell})\}_{\ell=0}^{m}$, for each model. Corrupted trajectories were then generated by randomly selecting a fraction of the snapshots and perturbing them with additive Gaussian noise. The magnitude of the perturbation was controlled by a noise level parameter $\eta \in (0,1)$. Specifically, given a pair
\[
(q,\eta)=(\text{outlier percentage},\,\text{noise level})\in[0,1]^2,
\]
we selected uniformly at random a subset
\[
\mathcal{I}\subset\{0,\ldots,m\}
\]
with cardinality \(|\mathcal{I}|=\lfloor q(m+1)\rfloor\). The corrupted trajectory
\[
D=\{\hat{\boldsymbol{x}}(t_{\ell})\}_{\ell=0}^{m}
\]
was then constructed according to rule
\[
\hat{x}_{i}(t_{\ell})=
\begin{cases}
x_{i}(t_{\ell})+\sigma_i\,\epsilon_{i,\ell}, & \text{if } \ell\in\mathcal{I},\\[4pt]
x_{i}(t_{\ell}), & \text{otherwise},
\end{cases}
\]
where {\(\epsilon_{i,\ell}\stackrel{\mathrm{i.i.d.}}{\sim}\mathcal{N}(0,1)\)}, and
\[
\sigma_i
=
\eta
\sqrt{
\frac{1}{m+1}
\sum_{\ell=0}^{m}
x_i(t_{\ell})^2
}.
\]
Thus, the noise amplitude for each state variable was scaled by its root-mean-square value along the clean trajectory \(\{\boldsymbol{x}(t_{\ell})\}_{\ell=0}^{m}\). Note that only the samples in $\mathcal{I}$ are perturbed, while the remaining samples are not corrupted by noise; hence, the noise level $\eta$ controls the magnitude of the outliers, and the outlier percentage $q$ controls their number.

We assess the performance of the SINDy variants as follows. For each model and each parameter pair \((q,\eta)\), we generate 100 independently corrupted trajectories $D=\{\hat{\boldsymbol{x}}(t_{\ell})\}_{\ell=0}^{m}$ using the procedure described above. Each trajectory is then used as input to the identification algorithm, yielding an estimated coefficient matrix ${\hat{\alpha}}$ in the chosen feature space. The recovered support is compared with the ground-truth support by classifying each candidate library term as either active (nonzero coefficient) or inactive (zero coefficient). For a single corrupted trajectory, the \textit{support recovery accuracy} is defined as
\begin{equation}
\label{eq:accuracy}
\mathrm{Acc}
=
\frac{\mathrm{TP}+\mathrm{TN}}
{\mathrm{TP}+\mathrm{TN}+\mathrm{FP}+\mathrm{FN}},
\end{equation}
where \(\mathrm{TP}\), \(\mathrm{TN}\), \(\mathrm{FP}\), and \(\mathrm{FN}\) denote the numbers of true positives, true negatives, false positives, and false negatives, respectively. Let \(\mathrm{Acc}_s(q,\eta)\) denote the accuracy obtained from the \(s\)-th trial, \(s=1,\ldots,100\), for a fixed parameter pair \((q,\eta)\). The \textit{mean support recovery accuracy} is defined as
\[
\overline{\mathrm{Acc}}(q,\eta)
=
\frac{1}{100}
\sum_{s=1}^{100}
\mathrm{Acc}_s(q,\eta).
\]
In addition, we quantify exact support recovery through the \textit{exact recovery rate}
\[
R_{\mathrm{exact}}(q,\eta)
=
\frac{1}{100}
\sum_{s=1}^{100}
\mathbf{1}_{\{\mathrm{Acc}_s(q,\eta)=1\}},
\]
which corresponds to the fraction of trials in which the sparsity pattern of the true coefficient vector is recovered exactly. 

Finally, to assess the quality of the recovered coefficients, we compute the relative coefficient error
\begin{equation}
\label{eq:coeff_error}
E_s(q,\eta) = \frac{\|\alpha- \hat{\alpha}_s\|_F}{\|\alpha\|_F},
\end{equation}
where $\alpha$ is the ground-truth coefficient matrix, $\hat{\alpha}_s$ is the matrix estimated in the $s$-th trial, and \mbox{$\|\cdot\|_F$} denotes the Frobenius norm. We report the median over the $100$ trials,
\[
\widetilde{E}(q,\eta) = \operatorname{median}\{E_1(q,\eta),\ldots,E_{100}(q,\eta)\}.
\]
Unlike the support-based metrics, this quantity is sensitive to the magnitude of the identified coefficients and therefore penalizes models that recover the correct sparsity pattern but with inaccurate values.

The mean support recovery accuracy and the exact recovery rate are evaluated for each model over all parameter pairs with
\begin{align*}
q,\eta \in  &\{2.5\%,\,5\%,\,7.5\%,\,10\%,\\
            &\,12.5\%,\,15\%,\,17.5\%,\,20\%\},
\end{align*}
yielding $64$ parameter combinations, and are presented in the form of heat maps to visualize their dependence on the outlier percentage and noise level. Since, for each pair $(q,\eta)$, we perform $100$ independent trials, our results are based on a total of $6,\!400$ corrupted trajectories per model, providing reliable estimates across all considered parameter pairs $(q,\eta)$. We compare the E-SINDy implementations provided by the PySINDy package\footnote{\url{https://github.com/dynamicslab/pysindy}.} with default parameters against our own Python implementation of SINDy and ILTS-SINDy\footnote{\url{https://github.com/fabioamaral08/ILTS-SINDy}.} using this synthetic dataset. All experiments were conducted {on a Ubuntu 22.04.4 LTS machine with an Intel(R) Xeon(R) CPU E5-2690 2.90GHz processor and 32 GB of RAM}.

\subsection{Hyperparameters}\label{sec:hyperparameters}

This section details the hyperparameter selection for the numerical experiments. All methods use STLS as the optimizer, whose only hyperparameter is the sparsity threshold $\lambda$. Since the ground-truth models of the benchmarks are known, $\lambda$ was tuned for each problem, and the selected values are reported in Table~\ref{tab:models_params}. Specifically, $\lambda$ was chosen as a round value no larger than half the magnitude of the smallest nonzero coefficient of the true model, and the same value was used for all compared methods, so that no method is favored by this choice.

Regarding the candidate mapping, all variants considered in our experiments employ a library
\[
\theta: \mathbb{R}^{N} \to \mathbb{R}^{d}, \qquad
\theta(\boldsymbol{z}) = \big(\theta_{1}(\boldsymbol{z}), \ldots, \theta_{d}(\boldsymbol{z})\big),
\]
whose components are the candidate functions used for system identification. Unless stated otherwise, we adopt the set of monomials of degree at most three. For the two-dimensional benchmarks (Lotka--Volterra and Van der Pol), this gives

\[
\theta(\boldsymbol{z}) =
(1, z_{1}, z_{2}, z_{1}^{2}, z_{1}z_{2}, z_{2}^{2}, z_{1}^{3}, z_{1}^{2}z_{2}, z_{1}z_{2}^{2}, z_{2}^{3}),
\]
whereas for the three-dimensional benchmarks (Lorenz and R\"ossler),

\begin{align*}
    \theta(\boldsymbol{z}) = \ &(1, z_{1}, z_{2}, z_{3}, z_{1}^{2}, z_{1}z_{2}, z_{1}z_{3}, z_{2}^{2}, z_{2}z_{3}, z_{3}^{2}, z_{1}^{3}, \\z_{1}^{2}z_{2},
    &\ z_{1}^{2}z_{3}, z_{1}z_{2}^{2}, z_{1}z_{2}z_{3}, z_{1}z_{3}^{2}, z_{2}^{3}, z_{2}^{2}z_{3}, z_{2}z_{3}^{2}, z_{3}^{3}).
\end{align*}

The SIR benchmark requires a more restrictive library because of the structure of the data. The state satisfies the constraint $z_{1} + z_{2} + z_{3} = 1$, which makes the columns of the full polynomial library linearly dependent: the constant term coincides with the sum of the linear terms, and each linear term $z_{i}$ equals $z_{i}^{2} + \sum_{j \neq i} z_{i}z_{j}$. To ensure that the candidate matrix has full column rank, we therefore discard the constant and the squared terms, retaining only
\[
\theta(\boldsymbol{z}) =
(z_{1}, z_{2}, z_{3}, z_{1}z_{2}, z_{1}z_{3}, z_{2}z_{3}).
\]

Finally, since the ABC flow is driven by trigonometric terms, polynomial candidates are inadequate. We instead use a trigonometric library with two frequencies,
\begin{align*}
\theta(\boldsymbol{z}) = \ &
(1, \cos z_1, \sin z_1, \cos z_2, \sin z_2,\\
&\ \cos z_3, \sin z_3, \cos 2z_1, \sin 2z_1,\\
&\ \cos 2z_2, \sin 2z_2, \cos 2z_3, \sin 2z_3).
\end{align*}

Turning to the method-specific parameters, preliminary experiments with ILTS-SINDy led us to set the target number of inliers in all tests to
\[
p = (m+1) - \left\lfloor 3\,q\,(m+1) \right\rfloor,
\]
where $q$ is the fraction of outliers and $m+1$ the number of samples. In the benchmark experiments, $q$ is the outlier fraction used to generate the corrupted data, so this choice assumes that the contamination level is known; the case in which it is unknown is considered in Section~\ref{sec:hyper_est}. The factor $3$ is not arbitrary: it follows from the use of central finite differences for derivative estimation, since a single corrupted observation contaminates up to three derivative estimates, namely the one at the affected time instant and those at its two neighbors. An ablation study supporting this choice is reported in Appendix~\ref{app:ablation}.

For E-SINDy and W-SINDy, we use the PySINDy implementations with their default settings. E-SINDy is run with 20 ensemble members, each fitted to 1001 samples drawn from the data with replacement (i.e., as many samples as the original dataset), and the individual models are aggregated by the coordinate-wise median (bragging). W-SINDy uses $K = 100$ test functions of the form given in Eq.~\eqref{eq:poly_test} with $r = s = 4$, so that each test function is a polynomial of degree $8$ on its support. Each support has half-width equal to $1/20$ of the observation window, i.e., a nominal width of $(T_f-T_0)/10$, and the centers are drawn uniformly at random over the region in which a full-width support fits inside the window. Since the placement of the supports is stochastic, the centers are resampled at every repetition, so the reported metrics also average over the test function placement.
\subsection{Numerical Results and Comparison}

\begin{figure}[htb!]
    \centering
    \includegraphics[width=\linewidth]{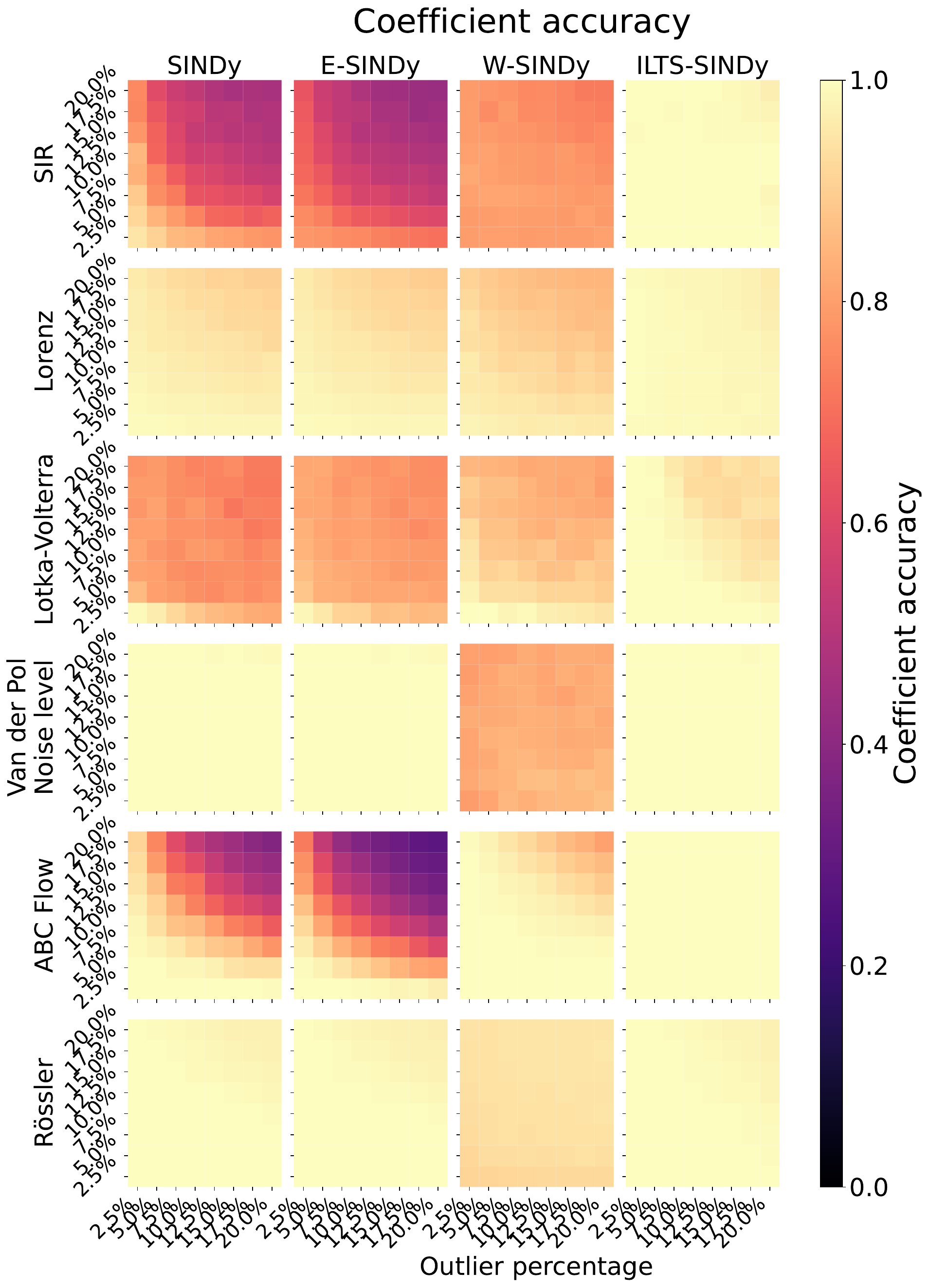}
    \caption{Mean support recovery accuracy $\overline{\mathrm{Acc}}(q,\eta)$ for different outlier fractions $q$ and noise levels $\eta$, for each method and benchmark system. Each cell reports the mean over $100$ corrupted realizations.}
    \label{fig:heatmap_all}
\end{figure}

Figure \ref{fig:heatmap_all} presents the accuracy heatmaps for SINDy, E-SINDy, W-SINDy and ILTS-SINDy across different combinations of noise level (outlier magnitude) $\eta$ and outlier percentage $q$ for the systems presented in Table \ref{tab:models_params}.

Overall, accuracy decreases as the level of data corruption increases, with the noise level $\eta$, i.e., the magnitude of the outliers, having a more pronounced effect than the outlier percentage $q$. Standard SINDy exhibits the greatest sensitivity to corruption, showing substantial performance degradation, particularly for the SIR, Lotka–Volterra and ABC Flow systems. W-SINDy demonstrates improved robustness and generally maintains moderate-to-high accuracy across the parameter space. However, ILTS-SINDy achieves the highest accuracy in almost all tested settings, for the systems and methods compared in this work, with values remaining close to 1 even under the most challenging settings. The only exceptions are a few settings of the R\"ossler system with noise levels up to $15\%$ and two low-noise settings of the Lorenz system, in which SINDy or E-SINDy attain a slightly higher accuracy (by at most $0.008$). The relatively uniform heatmaps observed for ILTS-SINDy indicate strong robustness to both noise level and outlier percentage, suggesting that the proposed method can reliably reconstruct system dynamics despite significant data corruption.

Table~\ref{tab:accuracy_noise_0.2} reports the numerical values underlying the top row of the colormap, corresponding to an outlier magnitude of $20\%$. ILTS-SINDy attains high accuracy with low variance across all benchmarks.

\begin{table*}
    \centering
    \caption{Support recovery accuracy $\mathrm{Acc}$ (mean $\pm$ standard deviation) at noise level $\eta = 20\%$, aggregated over all outlier fractions $q$. Values are means over the $100$ realizations of each configuration. The best result for each system is shown in bold.}

        \begin{tabular}{lcccc}
        \toprule
        Problem & SINDy & E-SINDy & W-SINDy & ILTS-SINDy \\
        \midrule
        SIR & $0.54 \pm 0.15$ & $0.49 \pm 0.11$ & $0.76 \pm 0.07$ & $\mathbf{0.99 \pm 0.03}$ \\
        Lorenz & $0.92 \pm 0.04$ & $0.92 \pm 0.03$ & $0.87 \pm 0.05$ & $\mathbf{0.98 \pm 0.02}$ \\
        Lotka-Volterra & $0.75 \pm 0.12$ & $0.78 \pm 0.07$ & $0.83 \pm 0.10$ & $\mathbf{0.95 \pm 0.07}$ \\
        Van der Pol & $0.997 \pm 0.02$ & $0.997 \pm 0.02$ & $0.81 \pm 0.09$ & $\mathbf{0.999 \pm 0.01}$ \\
        ABC Flow & $0.56 \pm 0.20$ & $0.41 \pm 0.16$ & $0.90 \pm 0.08$ & $\mathbf{1.00 \pm 0.00}$ \\
        R\"ossler & $0.98 \pm 0.02$ & $0.98 \pm 0.01$ & $0.95 \pm 0.02$ & $\mathbf{0.99 \pm 0.02}$ \\
        \bottomrule
    \end{tabular}

    \label{tab:accuracy_noise_0.2}
\end{table*}

\begin{figure}[htb!]
    \centering
    \includegraphics[width=\linewidth]{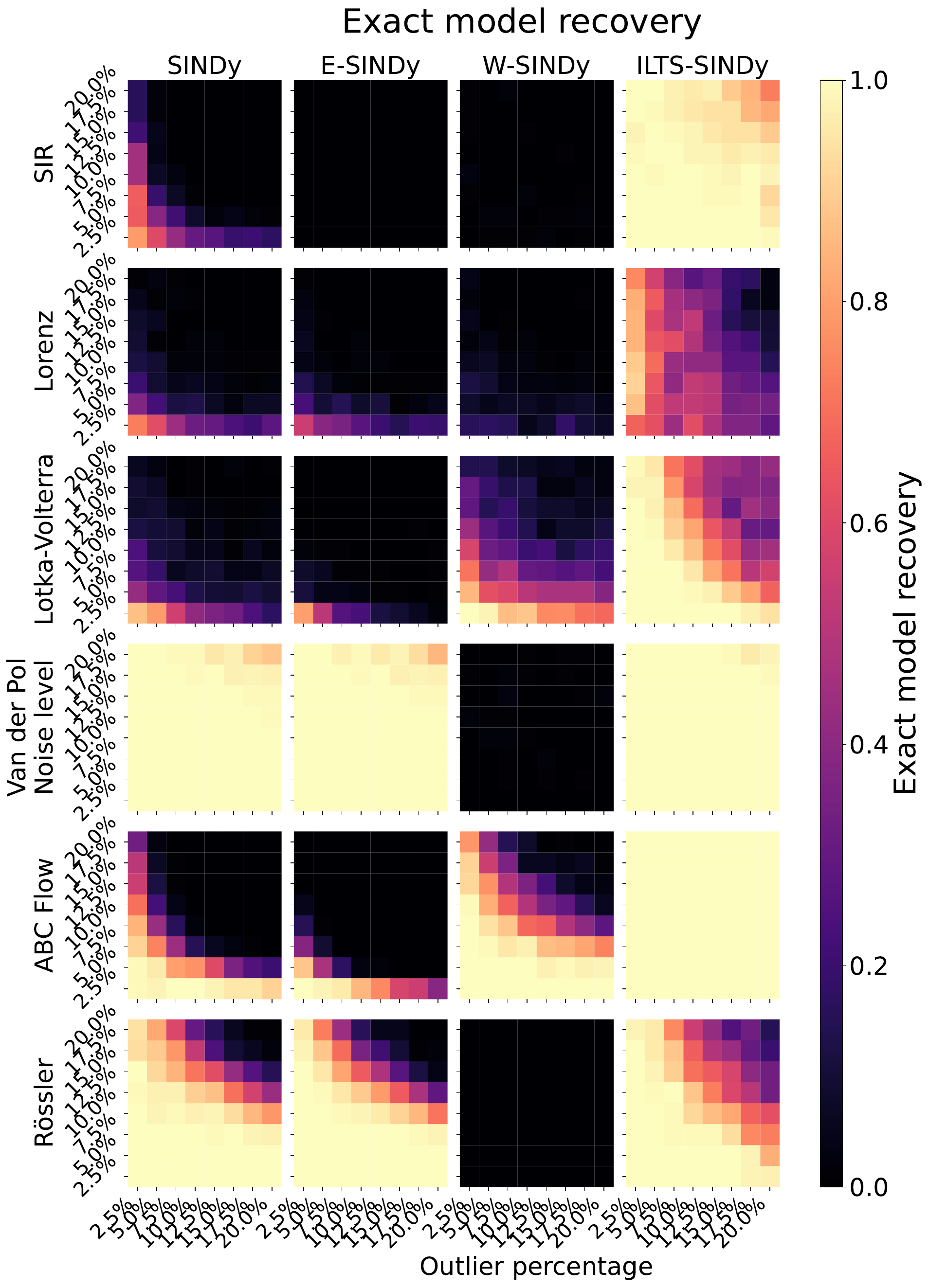}
    \caption{Exact recovery rate $R_{\mathrm{exact}}(q,\eta)$ for different outlier fractions $q$ and noise levels $\eta$, for each method and benchmark system. Each cell reports the fraction of the $100$ corrupted realizations in which the sparsity pattern of the true coefficient matrix was recovered exactly.}
    \label{fig:heatmap_hits_all}
\end{figure}

Figure \ref{fig:heatmap_hits_all} shows the exact recovery heatmaps, which quantify the percentage of trials in which the exact governing equations are recovered. 

In contrast to the accuracy results, the exact recovery heatmaps reveal substantial differences in the ability of the methods to recover the true model structure. This metric is deliberately strict: it takes the value $0$ whenever a single spurious coefficient is selected or a single true coefficient is missed, so values necessarily degrade in the corner corresponding to high noise and a high fraction of outliers. Standard SINDy achieves moderate identification rates only under low-noise and low-outlier conditions, deteriorating rapidly as corruption increases. Although E-SINDy frequently attains high mean accuracy (Figure~\ref{fig:heatmap_all}), its exact identification rates remain low across most settings, indicating that accurate trajectory reconstruction does not necessarily correspond to successful model discovery. The same dissociation is observed for W-SINDy, most notably on the Van der Pol oscillator and Rössler system, where competitive accuracy coexists with poor structural recovery. ILTS-SINDy consistently achieves the highest exact recovery rates across the six benchmarks, retaining the exact sparsity pattern under severe corruption for the SIR, Van der Pol, and ABC flow systems. The ABC flow is particularly noteworthy: Table~\ref{tab:exact_recovery_noise_0.2} reports a recovery rate of $100\%$ with zero standard deviation, meaning that the exact model was identified in every noise realization. These results indicate that ILTS-SINDy not only achieves a high support recovery accuracy but also more reliably recovers the underlying governing equations.

The exact recovery rate is a very strict, binary metric: a model is counted as a success only if its sparsity pattern exactly matches that of the true system. A single spurious term, or a single missing one, is enough for the model to be counted as a failure, regardless of how accurate the remaining coefficients are.

In the presence of outliers, standard least-squares fitting is strongly affected by the contaminated samples, which bias the estimated coefficients. As a result, some coefficients that should be zero remain above the sparsity threshold $\lambda$, and spurious terms appear in the identified model. Since ILTS-SINDy discards the samples with the largest residuals before fitting the final model, the outliers have much less influence on the regression, and spurious terms are far less likely to survive the thresholding step. For E-SINDy, the bootstrap subsamples contain, on average, the same fraction of outliers as the full dataset, so the aggregation over multiple models does not remove their influence.

This also explains why the difference between the methods is smaller in the accuracy heatmaps (Figure~\ref{fig:heatmap_all}) than in the exact recovery heatmaps (Figure~\ref{fig:heatmap_hits_all}): a model with a few small spurious coefficients may still achieve a reasonable accuracy, but it fails the exact recovery criterion. In other words, the remaining noise mainly affects the precision of the estimated coefficients of ILTS-SINDy, while its effect on the identified sparsity pattern is much smaller.

\begin{table*}
    \centering
    \caption{Exact recovery rate $R_{\mathrm{exact}}$ (mean $\pm$ standard deviation) at noise level $\eta = 20\%$, aggregated over all outlier fractions $q$. Values are means over the $100$ realizations of each configuration. The best result for each system is shown in bold.}
    
    \begin{tabular}{lcccc}
        \toprule
        Problem & SINDy & E-SINDy & W-SINDy & ILTS-SINDy \\
        \midrule
        SIR & $0.02 \pm 0.15$ & $0.00 \pm 0.00$ & $0.004 \pm 0.06$ & $\mathbf{0.92 \pm 0.27}$ \\
        Lorenz & $0.01 \pm 0.08$ & $0.001 \pm 0.04$ & $0.01 \pm 0.07$ & $\mathbf{0.34 \pm 0.47}$ \\
        Lotka-Volterra & $0.02 \pm 0.13$ & $0.00 \pm 0.00$ & $0.07 \pm 0.26$ & $\mathbf{0.62 \pm 0.49}$ \\
        Van der Pol & $0.96 \pm 0.19$ & $0.96 \pm 0.20$ & $0.004 \pm 0.06$ & $\mathbf{0.99 \pm 0.09}$ \\
        ABC Flow & $0.04 \pm 0.21$ & $0.00 \pm 0.00$ & $0.18 \pm 0.38$ & $\mathbf{1.00 \pm 0.00}$ \\
        R\"ossler & $0.36 \pm 0.48$ & $0.30 \pm 0.46$ & $0.00 \pm 0.00$ & $\mathbf{0.55 \pm 0.50}$ \\
        \bottomrule
    \end{tabular}
    
    \label{tab:exact_recovery_noise_0.2}
\end{table*}

\begin{figure}
    \centering
    \includegraphics[width=\linewidth]{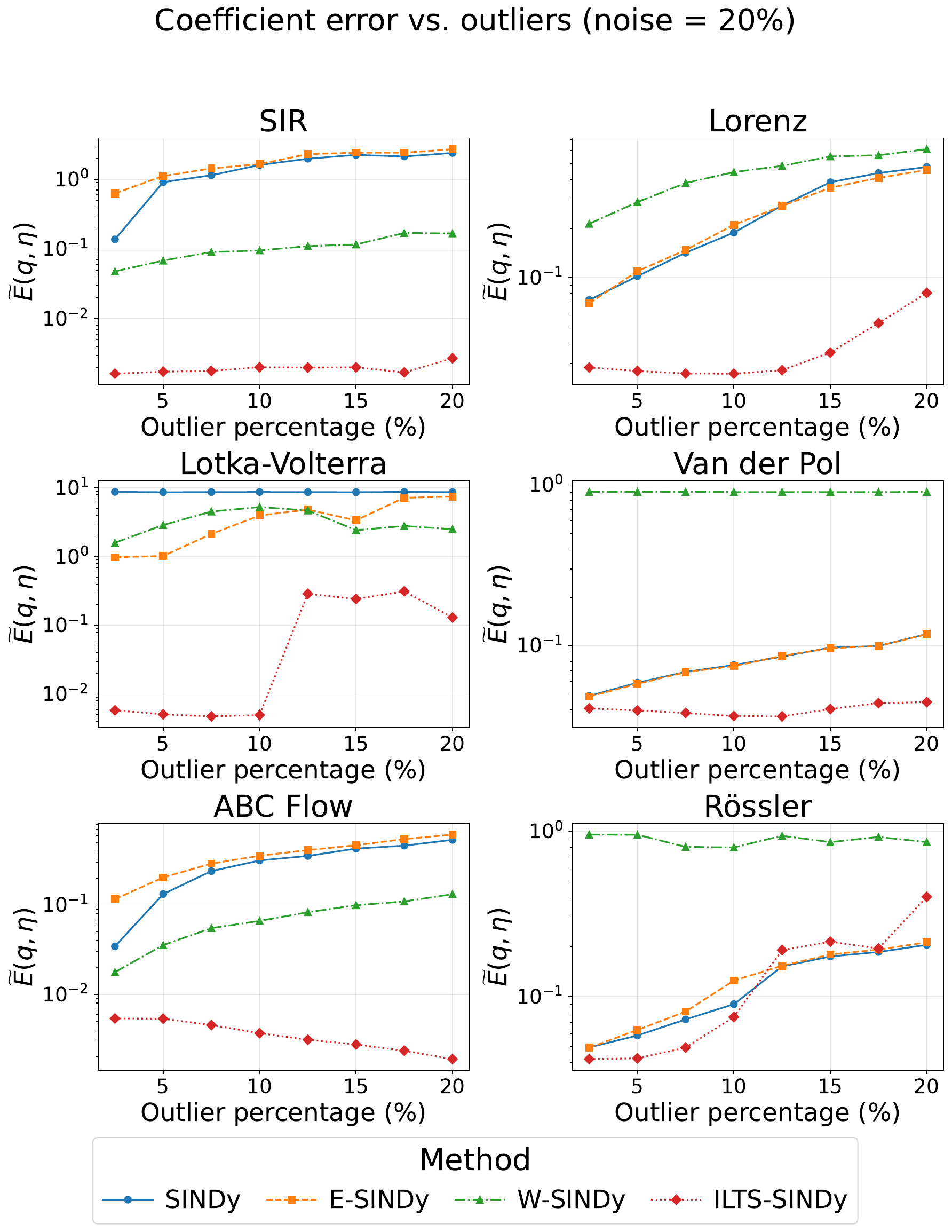}
    \caption{Median relative coefficient error $\widetilde{E}(q,\eta)$ as a function of the outlier fraction $q$, at the highest noise level tested ($\eta = 20\%$), for the six benchmark systems. Each point is the median over $100$ corrupted realizations.}
    \label{fig:coeff_err}
\end{figure}

Turning to the coefficients themselves, Figure~\ref{fig:coeff_err} shows the median coefficient error $\widetilde{E}(q,\eta)$ as a function of the outlier percentage, at the highest noise level tested. ILTS-SINDy attains the lowest error on five of the six benchmarks, by up to two orders of magnitude on the SIR and ABC flow systems, and its error is nearly insensitive to the outlier percentage for the SIR and Van der Pol systems. For the Lotka--Volterra system, the error increases at $q=12.5\%$, but remains at least one order of magnitude below that of the competing methods. The exception is the R\"ossler system, where ILTS-SINDy is the most accurate method up to approximately $10\%$ outliers but is overtaken beyond this point, ending at roughly twice the error of SINDy and E-SINDy at $20\%$; all three nevertheless remain well below W-SINDy over the entire range. Thus, although the performance of ILTS-SINDy may deteriorate under heavier contamination for certain systems, the overall results remain favorable to the proposed method. The small non-monotone variations observed in some curves, such as those of the Lotka--Volterra system for $q\ge12.5\%$ and of the R\"ossler system between $12.5\%$ and $17.5\%$, can be attributed to the randomness of the corrupted samples in each realization, since the corresponding values are close to each other.

We hypothesize that the unexpected behavior observed for the ABC flow, whose error remains low and even decreases slightly as the outlier percentage increases, is related to the structure of its library, which differs from those used by the other systems. Since the trigonometric candidate functions are bounded, outliers cause smaller changes in the corresponding rows of the library matrix, which may explain the higher accuracy of the estimated coefficients. By contrast, in a polynomial library, outliers can have a larger effect due to the higher-order terms.

Finally, we emphasize that the main goal of ILTS-SINDy, as of SINDy in general, is the identification of the structure of the model, measured by the accuracy and exact recovery rate of Figures~\ref{fig:heatmap_all} and~\ref{fig:heatmap_hits_all}. The coefficient values are a secondary estimate: since the trimming does not guarantee that all corrupted samples are removed, the coefficients may still be affected by the remaining contamination even when the structure is correctly identified. Once the structure is known, the coefficients can be re-estimated in the reduced function space.

\begin{figure}
    \centering
    \includegraphics[width=\linewidth]{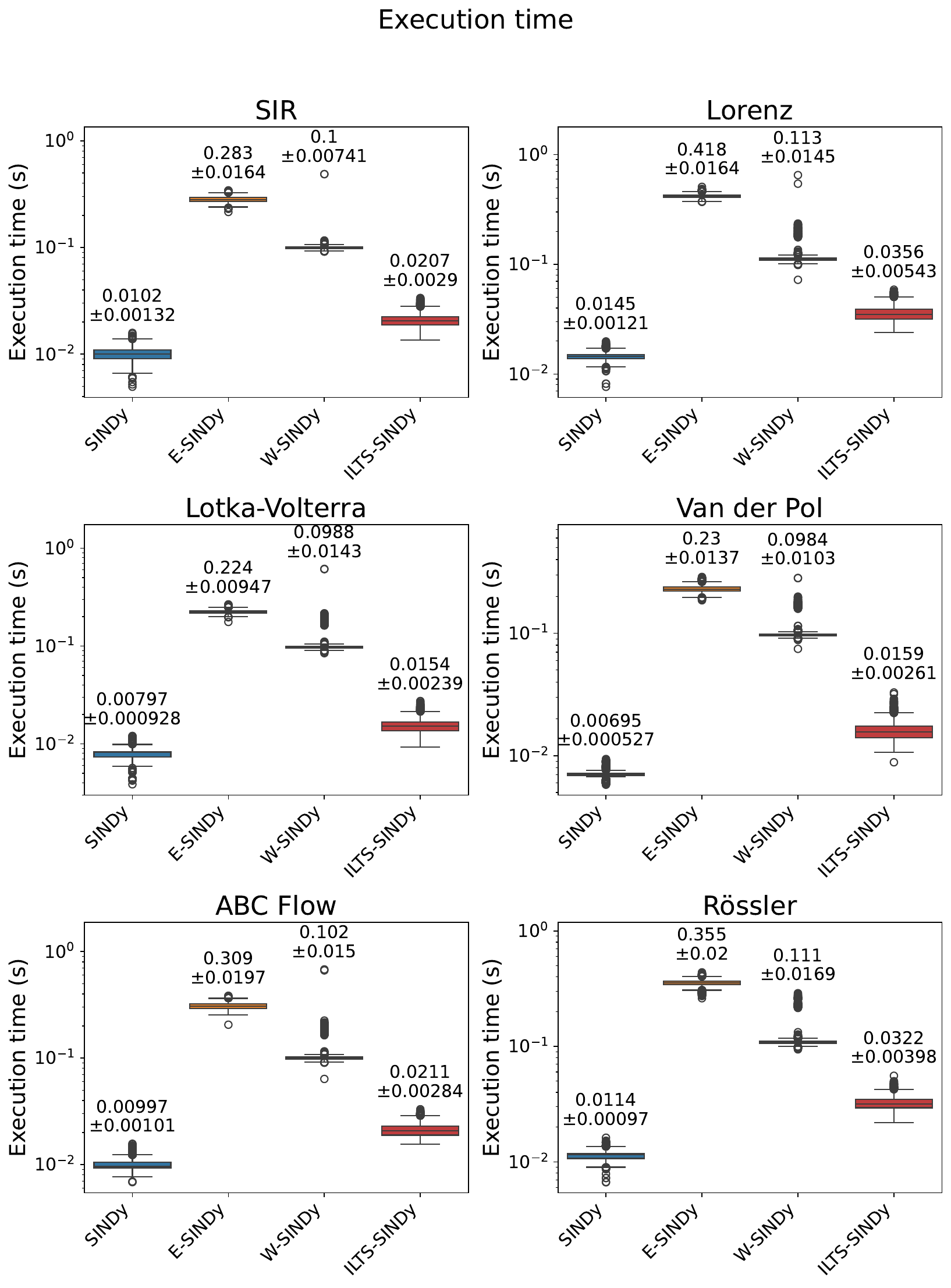}
    \caption{Execution time of each method across the six benchmark systems. Box plots aggregate all $6400$ runs per system, obtained by pooling every combination of outlier fraction $q$ and noise level $\eta$. Annotations  report the mean and standard deviation in seconds}
    \label{fig:time_benchmark}
\end{figure}

The final benchmark concerns the computational cost of each method. Since the runtime is essentially unaffected by the noise level and the outlier fraction, we aggregate the statistics over all $6400$ simulations of each problem, yielding the box plots of Figure~\ref{fig:time_benchmark}, where the annotations report the mean execution time and its standard deviation. As expected, standard SINDy is the fastest method, since it constitutes the baseline on which the remaining variants are built, with mean times between $7.0$~ms (Van der Pol) and $14.5$~ms (Lorenz). ILTS-SINDy ranks second, from $15.4$~ms to $35.6$~ms, corresponding to an overhead of roughly two to three times the cost of standard SINDy. W-SINDy requires approximately $100$~ms on every benchmark ($98.4$--$113$~ms), three to six times more than ILTS-SINDy, and E-SINDy is the most expensive, ranging from $224$~ms to $418$~ms, since it fits one model per ensemble member at each simulation. The overhead introduced by ILTS-SINDy is therefore modest in absolute terms, and it remains an order of magnitude faster than E-SINDy while delivering substantially better accuracy and support recovery.

\subsection{Hyperparameter estimation}\label{sec:hyper_est}

Algorithm~\ref{alg:ilts-sindy} takes as input the number of inliers $p$ and the sparsity threshold $\lambda$, the latter being the threshold used by STLS to promote sparsity. In practice, applying Algorithm~\ref{alg:ilts-sindy} requires estimating $p$ and tuning $\lambda$, neither of which is known a priori. In this section, we describe automatic tuning procedures for $p$ and $\lambda$, and illustrate their performance on the SIR model.

For the estimation of $p$, we can use the following heuristic suggested in \cite{Andreani2}. Given a user-defined interval $[p_{\min}, p_{\max}]$ of candidate values, with $p_{\min} < p_{\max}$, we select
\begin{equation}
\label{eq:p_selection}
p^{*} = \operatorname*{arg\,max}_{p \in [p_{\min},\, p_{\max}]} \frac{L(p+1) - L(p)}{L(p)},
\end{equation}
where
\[
L(p) = \tfrac{1}{2} \big\| \hat{A}\hat{\alpha}_i - \hat{b}_i \big\|_2^2,\ \]
and
\[
(\hat{\alpha}_i, \hat{A}_i, \hat{b}_i) = \mathrm{ILTS}(A, b_i, p).
\]
The rationale is as follows. As long as $p$ remains below the true number of inliers, the observations added to the inlier set when $p$ is incremented are themselves consistent with the underlying model, so the objective $L(p)$ increases only mildly. Once $p$ exceeds the true number of inliers, however, further increments force corrupted observations into the set, producing a pronounced increase in $L(p)$. The criterion in Eq.~\eqref{eq:p_selection} therefore identifies the point at which the growth of $L(p)$ changes most abruptly, which serves as an estimate of the number of inliers.

\begin{figure}
    \centering
    \includegraphics[width=\linewidth]{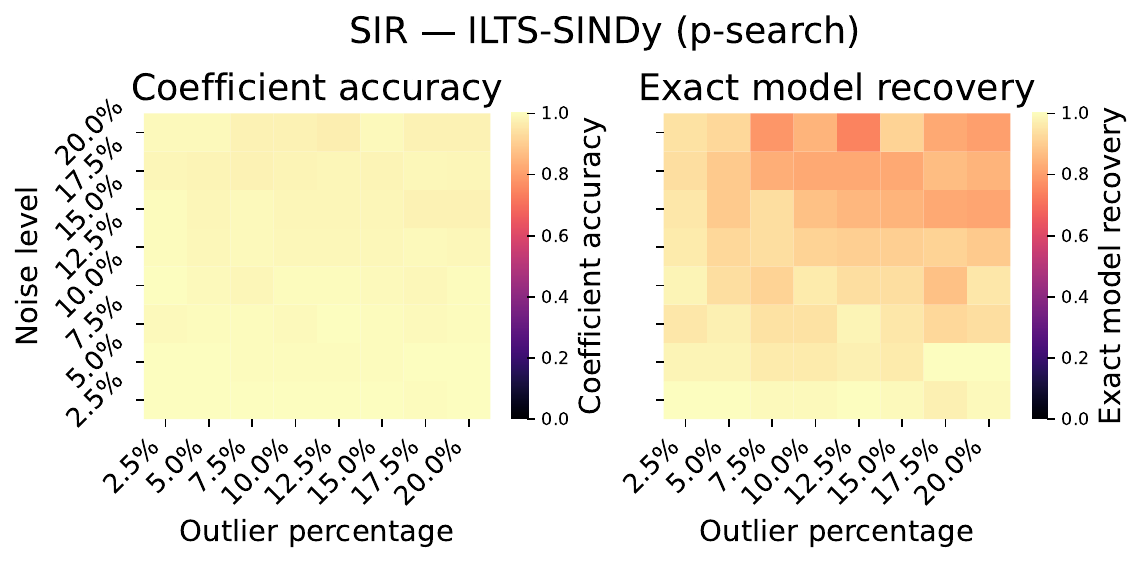}
    \caption{Mean support recovery accuracy (left) and exact recovery rate (right) for ILTS-SINDy on the SIR model with $p$ selected automatically by Eq.~\eqref{eq:p_selection} and $\lambda$ fixed as in Table~\ref{tab:models_params}. Each cell is a mean over $100$ realizations.}
    \label{fig:grid_search_p}
\end{figure}

Figure~\ref{fig:grid_search_p} reports the results obtained with the automatic selection of $p$ for the SIR model. The candidate values are constructed from an assumed inlier fraction ranging over $80\%$ to $100\%$ of the data: for each candidate count $\tilde{p} \in \{\lceil 0.8(m+1) \rceil, \ldots, m+1\}$, the corresponding target number of inliers is
\[
p = 3\tilde{p} - 2(m+1),
\]
and only positive values are retained. This is exactly the rule $p = (m+1) - \lfloor 3q(m+1) \rfloor$ written in terms of $\tilde{p} = (1-q)(m+1)$, so the factor $3$ again accounts for the propagation of a single corrupted sample to three central finite-difference estimates. Consequently, the candidate values of $p$ are spaced by $3$ and span the range from $0.4(m+1)$ to $m+1$. The resulting accuracy and exact recovery degrade only slightly with respect to the case in which the true value of $p$ is known, and remain superior to those of the competing methods reported in Figures~\ref{fig:heatmap_all} and~\ref{fig:heatmap_hits_all} for this system.

For the sparsity threshold, a similar criterion can be used. For a fixed inlier count $p_i^{*}$, let $J_i(\lambda)$ denote the residual of the STLS fit obtained with threshold $\lambda$ on the retained inlier set. Given a grid of candidate values $\lambda_1 < \cdots < \lambda_N$, logarithmically spaced so as to cover several orders of magnitude, we select
\begin{equation}
    \label{eq:lam_selection}
    k_i^{*} = \operatorname*{arg\,max}_{k \in \{1,\ldots,N-1\}} \frac{J_i(\lambda_{k+1}) - J_i(\lambda_k)}{\lambda_{k+1} - \lambda_k},
    \quad \lambda_i^{*} = \lambda_{k_i^{*}},
\end{equation}
that is, the difference quotient approximating the derivative of the residual with respect to the threshold. Increasing $\lambda$ progressively removes terms from the model. While only spurious terms are discarded, the residual grows slowly; once a genuine term is eliminated, it increases sharply. The criterion therefore locates this transition, and $\lambda_i^{*}$ is taken as the last threshold before it.

\begin{figure}
    \centering
    \includegraphics[width=\linewidth]{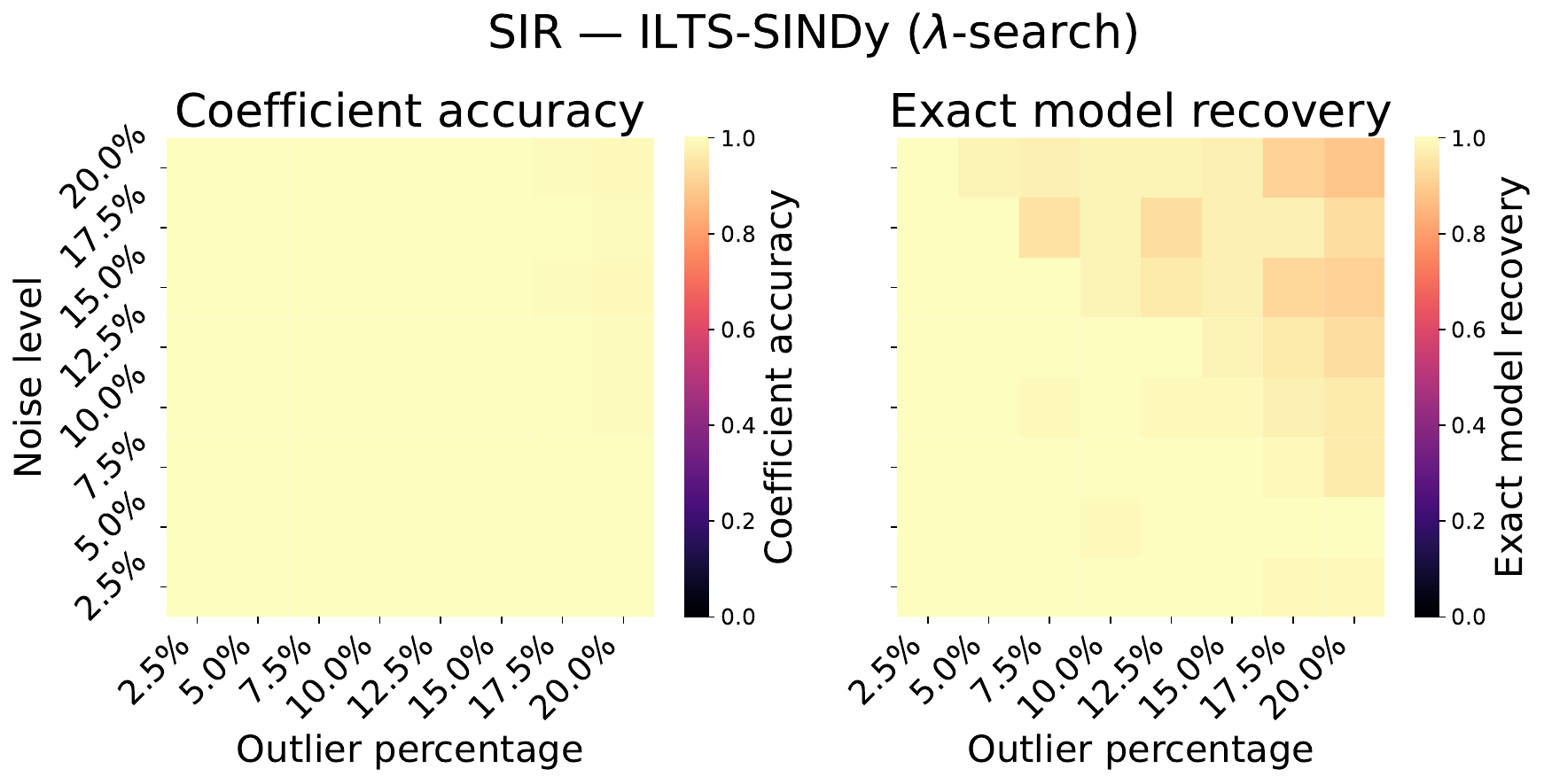}
    \caption{Mean support recovery accuracy (left) and exact recovery rate (right) for ILTS-SINDy on the SIR model with $\lambda$ selected automatically by Eq.~\eqref{eq:lam_selection} over a logarithmic grid in $[10^{-5}, 1]$, and $p$ assumed known. Each cell is a mean over $100$ realizations.}
    \label{fig:grid_search_lambda}
\end{figure}

Figure~\ref{fig:grid_search_lambda} reports the results of a grid search over $\lambda \in \{10^{-5}, 10^{-4}, \ldots, 10^{0}\}$ (one value per decade), performed independently for each state, for the SIR model, with $p$ assumed known. The automatically selected thresholds slightly improve upon the fixed value of $\lambda$ used in the previous experiments.

Both selection rules are applied independently for each state component $i$, so that the inlier count $p_i^{*}$ and the threshold $\lambda_i^{*}$ may differ across the equations of the system. This is deliberate, and it explains the improvement observed above: the components of a dynamical system generally differ in the magnitude of their coefficients and in the number of active library terms, so a single global threshold necessarily compromises between them, whereas a per-component threshold adapts the sparsity level to each equation.

\subsection{Application to real economic data}
\label{sec:realdata}
 
\begin{figure}[htbp]
  \centering
  \includegraphics[width=\linewidth]{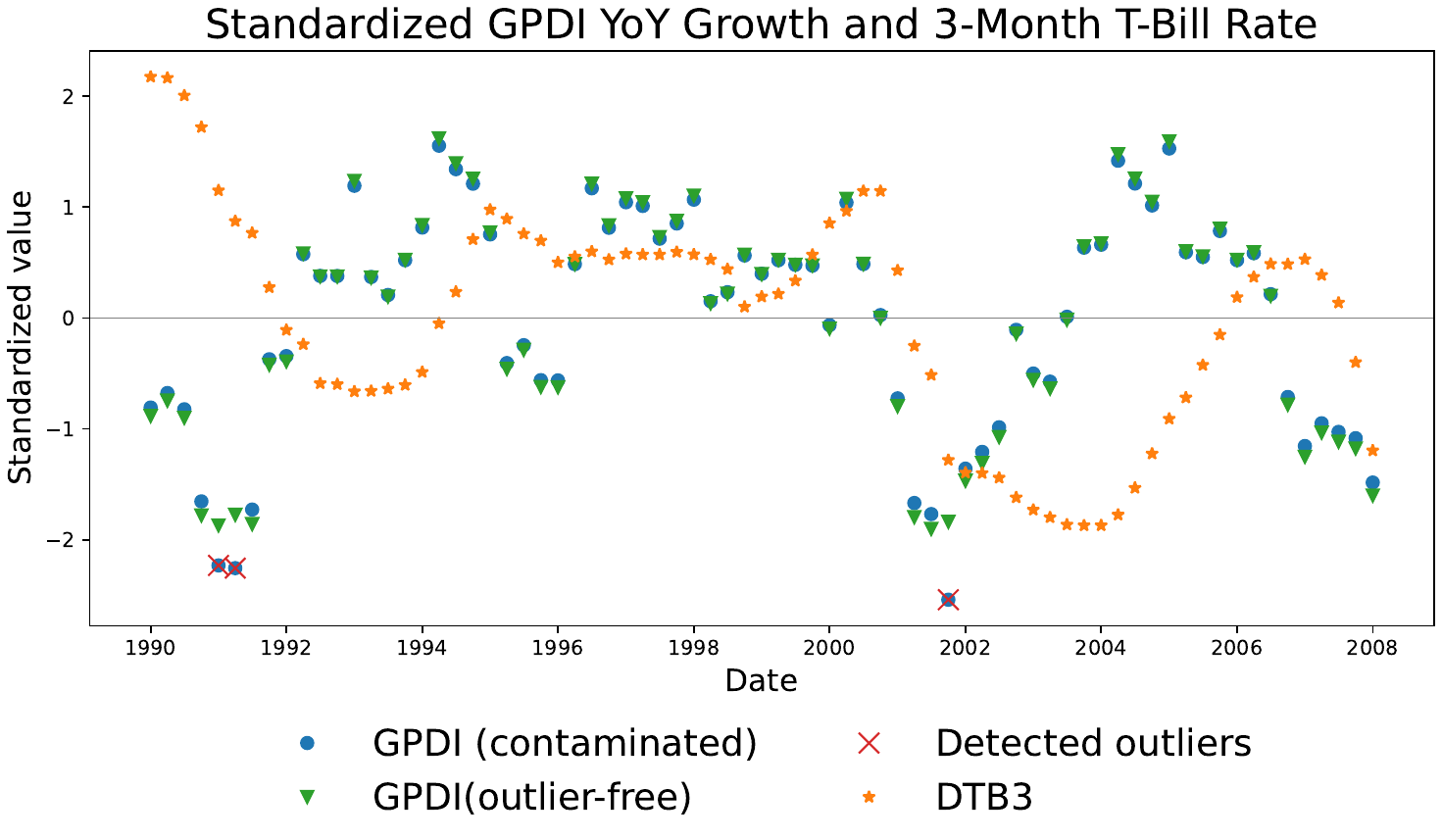}
  \caption{Standardized GPDI year-over-year growth and 3-month Treasury bill
  rate, 1990-Q1 to 2008-Q1: raw series, clean series (outliers replaced and
  Savitzky--Golay smoothed), and the three $z$-score outliers.}
  \label{fig:realdata_data}
\end{figure}

As an illustration, we apply ILTS-SINDy to the investment--interest-rate system studied in~\cite{li2025sparse}. The state is $x$, the year-over-year growth of US Gross Private Domestic Investment (Federal Reserve Economic Data, FRED: GPDI), and $y$, the 3-month Treasury bill rate (FRED: DTB3) averaged to quarters, both standardized. The period is 1990-Q1 to 2008-Q1 (73 samples, the points plotted in Fig.~7 of~\cite{li2025sparse}). The reference study replaced outliers detected by a $z$-score rule (three samples, in the 1991 and 2001 recessions, in our reproduction), smoothed the series with a Savitzky--Golay filter, and identified its model with E-SINDy. We therefore use two datasets: \emph{raw}, with the outliers kept and no smoothing, and \emph{clean}, preprocessed as in the reference (Fig.~\ref{fig:realdata_data}). Both use the library $[1, x, y, x^2, xy, y^2]$, time in years, and finite-difference derivatives; the reference coefficients, estimated at quarterly frequency, are rescaled to years. Since neither the outlier fraction nor the true model is known for these data, both tuning parameters of ILTS-SINDy were selected from the data, independently for each equation, using the criteria of Section~\ref{sec:hyper_est}: $p$ was selected by Eq.~\eqref{eq:p_selection}, with candidate values constructed as in that section from inlier fractions between $90\%$ and $99\%$, which retained $55$ and $70$ of the $73$ samples for the two equations on the raw data and $67$ and $61$ on the clean data; and $\lambda$ was selected by Eq.~\eqref{eq:lam_selection}, which gave $\lambda = 0.1$ in all four fits. For comparison, standard SINDy was run with threshold $0.001$, and E-SINDy with an ensemble of $100$ bagged STLS models.

The reference study states that its goal is ``not to perform precise forecasting but to uncover interpretable dynamical structures''~\cite[Sec.~4]{li2025sparse}. We keep the comparison in that scope: we compare the structure of the models and how closely their free simulation, started from the first sample, reproduces the observed series over the 18 years. Errors, measured by the root-mean-square error (RMSE) and the mean absolute error (MAE), are computed on all 73 samples, including those trimmed by ILTS-SINDy. 

\begin{table*}[htbp]
  \centering
  \caption{Identified models (per year; $x$: GPDI growth, $y$: Treasury bill
  rate, both standardized) and RMSE of the 18-year free simulation against the
  data. Best RMSE per dataset in bold.}
  \label{tab:realdata_models}
  \small
  \setlength{\tabcolsep}{4pt}
  
  \begin{tabular}{@{}llcccc@{}}
    \toprule
    Model & Equations & Terms & Cycle period & \multicolumn{2}{c}{RMSE} \\
    \cmidrule(l){5-6}
    (identified on) & & & (yr) & clean & raw \\
    \midrule
    \begin{tabular}[c]{@{}l@{}}Reference~\cite{li2025sparse}\\(clean)\end{tabular} &
      $\begin{cases}
        \dot{x} &= -0.528\,y - 0.048\,x^2 - 0.032\,xy \\
        \dot{y} &= -0.276 + 0.772\,x - 0.052\,xy + 0.052\,y^2
      \end{cases}$ & 7 & 9.7 & 1.327 & \textbf{1.339} \\[12pt]
    \begin{tabular}[c]{@{}l@{}}ILTS-SINDy\\(clean)\end{tabular} &
      $\begin{cases}
        \dot{x} &= -0.220 - 0.637\,y \\
        \dot{y} &= -0.573 + 0.939\,x + 0.221\,x^2 + 0.102\,xy + 0.168\,y^2
      \end{cases}$ & 7 & 7.4 & \textbf{1.116} & -- \\[12pt]
    \begin{tabular}[c]{@{}l@{}}ILTS-SINDy\\(raw)\end{tabular} &
      $\begin{cases}
        \dot{x} &= 0.220 - 0.694\,y - 0.142\,x^2 + 0.254\,xy \\
        \dot{y} &= -0.118 + 0.650\,x
      \end{cases}$ & 6 & 9.7 & -- & 1.730 \\
    \bottomrule
  \end{tabular}
\end{table*}
 
All models share the same structure: investment growth is driven down by the interest rate ($-y$ in $\dot{x}$), and the rate is driven up by investment growth ($+x$ in $\dot{y}$). This coupling produces a weakly damped oscillation (weakly amplified for raw ILTS-SINDy) around a single equilibrium, an endogenous business cycle, which the reference interprets as Lotka--Volterra-type feedback.
 
\begin{figure}[htbp]
  \centering
  \includegraphics[width=\linewidth]{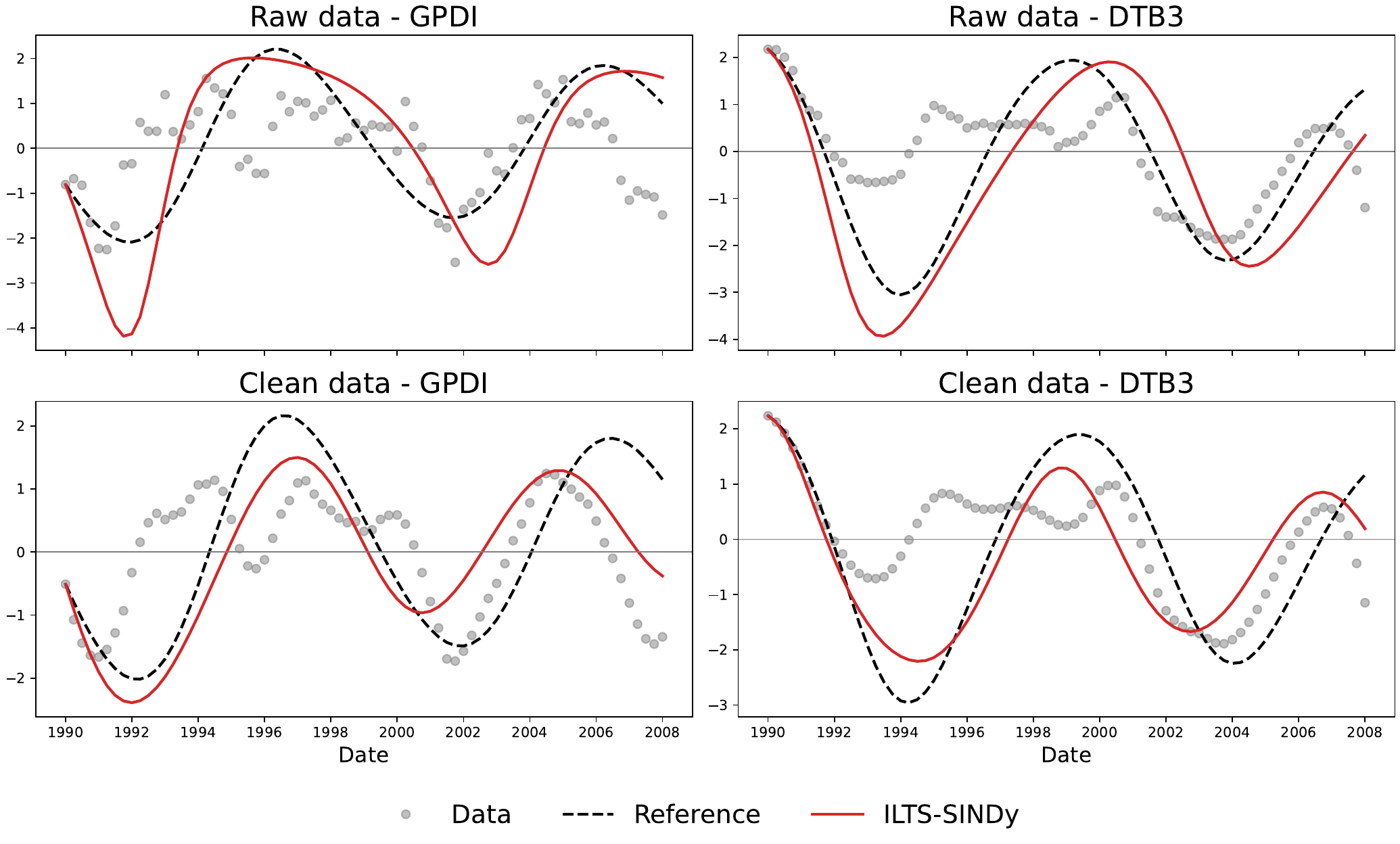}
  \caption{Free simulation over 18 years from the first sample: data,
  reference model~\cite{li2025sparse} and ILTS-SINDy, on the raw (top) and
  clean (bottom) datasets.}
  \label{fig:realdata_freerun}
\end{figure}

On the clean data, ILTS-SINDy fits the observed series better than the reference with the same number of terms: the RMSE falls by 16\% (1.116 vs 1.327) and the MAE from 1.031 to 0.851 (Fig.~\ref{fig:realdata_freerun}).

On the raw data, the fit of ILTS-SINDy is worse than that of the reference (1.730 vs 1.339). The comparison is not like for like, however: the reference was identified on cleaned and smoothed data, while ILTS-SINDy works directly on outlier-contaminated data and noisy derivatives. Under these conditions, standard SINDy and E-SINDy (the method used for the reference model) fit spurious quadratic terms and diverge in finite time (at $t = 13$ and 11 years). ILTS-SINDy stays bounded, gives the sparsest model (6 terms), keeps the same mechanism and recovers the reference cycle period almost exactly (9.68 vs 9.73 years). The samples it trims lie adjacent to each of the three outliers; thus, it identifies corrupted derivatives without requiring manual outlier rules. 

In practice, ILTS-SINDy therefore offers an acceptable, interpretable model straight from contaminated data, and a better-fitting one when the data can be cleaned.

\section{Conclusion}\label{sec:4}

In this work we proposed ILTS-SINDy, a filter-then-sparsify pipeline for the sparse identification of ODE models from possibly noisy trajectories. The method consists of a simple two-step procedure. First, we remove snapshots of the trajectory identified as outliers using the ILTS algorithm. Then, on the resulting cleaned trajectory, we seek a sparse approximate solution of the associated regression problem using the STLS algorithm. Decoupling outlier filtering from sparse regression improves the identification of the model structure from corrupted data. Under severe outlier contamination, it achieves, in most of the tested cases, exact structural identification and up to two orders of magnitude lower coefficient errors, while remaining an order of magnitude faster than ensemble methods. Consequently, ILTS-SINDy provides a fast and robust approach for automated model discovery in outlier-prone environments.

The present study also has limitations. First, in the benchmark experiments only the samples selected as outliers are perturbed, while the remaining samples are noise-free; the combination of outliers with measurement noise affecting all samples was not considered. Moreover, as is standard practice in SINDy benchmark studies, the tuning parameters $p$ and $\lambda$ were set using knowledge of the contamination level and of the true model. Second, the automatic parameter selection rules given by \eqref{eq:p_selection} and \eqref{eq:lam_selection} are heuristic procedures and are not guaranteed to identify the correct number of inliers or the optimal sparsity level, respectively. Third, although the structure of the model is identified in most cases, estimating accurate coefficient values is harder for some systems under heavy contamination, as shown in Section~\ref{sec:num_results} for the Lotka--Volterra and R\"ossler systems; the focus of the method remains the identification of the structure of the system. Finally, as in standard SINDy, the choice of the candidate library requires domain knowledge of the system under study, a limitation that ILTS-SINDy does not overcome.

While the proposed ILTS algorithm successfully automates outlier detection, the current framework identifies these anomalies only after numerical differentiation. This can be an issue, as a single outlier can propagate errors across the numerical stencil, inadvertently corrupting adjacent derivative estimates. A primary avenue for future work is the development of pre-differentiation outlier identification methods; by removing corrupted observations before computing derivatives, a significantly larger portion of the valid dataset can be preserved for model discovery. To achieve this and further enhance the current ILTS architecture, future efforts will explore integrating the LOVO framework \cite{Andreani2,Castelani}. As a broader generalization of LTS that accommodates non-linear models, LOVO offers a robust foundation for more complex outlier rejection tasks. Extending the SINDy pipeline with LOVO could provide a powerful mechanism for addressing pre-differentiation filtering and broadening the method's applicability to highly nonlinear dynamics. Lastly, establishing rigorous theoretical guarantees for the ILTS-SINDy estimators, such as consistency and support recovery, remains an open question and is left for future research.

\section*{Funding}
F. Amaral was supported by the São Paulo Research Foundation (FAPESP) process numbers 2021/07034-4 and 2023/06035-2. C. Oishi was partially supported by the São Paulo Research Foundation (FAPESP) process numbers  2024/02658-8 and 2023/14427-8, and the National Council for Scientific and Technological Development (CNPq), grant  307228/2023-1. G.N. Grapiglia was supported by FRS-FNRS, Belgium (Grant PINT-BILAT-M 40026198).

\section*{Acknowledgements}
The support and resources from the LSNIA-FCT/UNESP (\url{https://lsnia-unesp.github.io/}) cluster were instrumental in carrying out the simulations.

\section*{Competing Interests}
The authors have no relevant financial or non-financial interests to disclose.

\section*{Author Contributions}
Geovani N. Grapiglia and Cassio M. Oishi contributed to the study conception and design. Code implementation, and data production were performed by Fabio Amaral. All authors contributed to the analysis of the results and the writing of the manuscript. All authors read and approved the final manuscript.

\section*{Data Availability}

All code and implementation details for this work are publicly available at \url{https://github.com/fabioamaral08/ILTS-SINDy}. The repository includes scripts for generating all figures presented in this paper and can be used to reproduce the numerical results. The data used for the results are available at \url{https://doi.org/10.5281/zenodo.21040456}.

\begin{appendices}

\section{Ablation of the ILTS Multiplier Factor}\label{app:ablation}

A core feature of the proposed methodology is the outlier detection capacity of the ILTS algorithm, which requires a prior estimation of the number of trusted inliers ($p$). However, because the downstream SINDy pipeline utilizes a numerical differentiation stencil spanning three points, any single outlier corrupts its immediate temporal neighbors. Consequently, the expected number of outliers must be systematically overestimated to ensure the algorithm effectively isolates the uncorrupted data

To assess the true anomaly detection performance of the ILTS layer under this overestimation regime, we propose using the recall metric, rather than accuracy, defined as
\begin{equation}
\label{eq:recall}
\text{Recall} = \frac{TP}{TP + FN}.
\end{equation}
In this context, recall serves as a more relevant performance indicator. Since the deliberate overestimation of the outlier percentage inherently inflates the false positive rate, standard accuracy would be artificially suppressed. Utilizing recall isolates the framework's true capability to capture genuine anomalies.

In Table \ref{tab:outlier_recall}, we report the mean recall and standard deviation across 6,400 noise realizations for each of the three benchmark dynamical systems, utilizing the same dataset described in Section \ref{sec:num_results}.  Within this framework, the target number of trusted inliers is explicitly defined as:
\begin{equation}
\label{eq:inlier_count}
p = (m+1) - \lfloor \rho q (m+1) \rfloor,
\end{equation}
where $q \in (0, 1)$ represents the true outlier percentage, $m+1$ is the total number of sampled data points, and $\rho \ge 1$ is the multiplier factor designed to compensate for the temporal spread of anomalies caused by the numerical differentiation stencil. The 'Exact' row entry in Table \ref{tab:outlier_recall} represent the ideal theoretical value of $\rho$ required for each specific system to perfectly match the total quantity of corrupted points (the original anomalies combined with their adjacent contaminated neighbors).

\begin{table*}[htb!]
\caption{Recall in outlier detection for ILTS-SINDy}
\label{tab:outlier_recall}
\centering
\begin{tabular}{c|ccc}
\hline
\begin{tabular}[c]{@{}c@{}}Multiplier\\ Factor \(\rho\)\end{tabular} & SIR & Lorenz & \begin{tabular}[c]{@{}c@{}}Lotka-\\ Volterra\end{tabular} \\ \hline
1 & $0.871 \pm 0.032$ & $0.884 \pm 0.013$ & $0.877 \pm 0.018$ \\
2 & $0.952 \pm 0.010$ & $0.931 \pm 0.010$ & $\mathbf{0.923 \pm 0.011}$ \\
3 & $\mathbf{0.961 \pm 0.013}$ & $\mathbf{0.939 \pm 0.010}$ & $0.918 \pm 0.017$ \\
4 & $0.954 \pm 0.020$ & $0.926 \pm 0.012$ & $0.898 \pm 0.025$ \\
Exact & $\mathbf{0.964 \pm 0.011}$ & $\mathbf{0.943 \pm 0.009}$ & $\mathbf{0.929 \pm 0.014}$ \\ \hline
\end{tabular}
\end{table*}

\begin{figure*}[htb!]
    \centering
    \includegraphics[width=\linewidth]{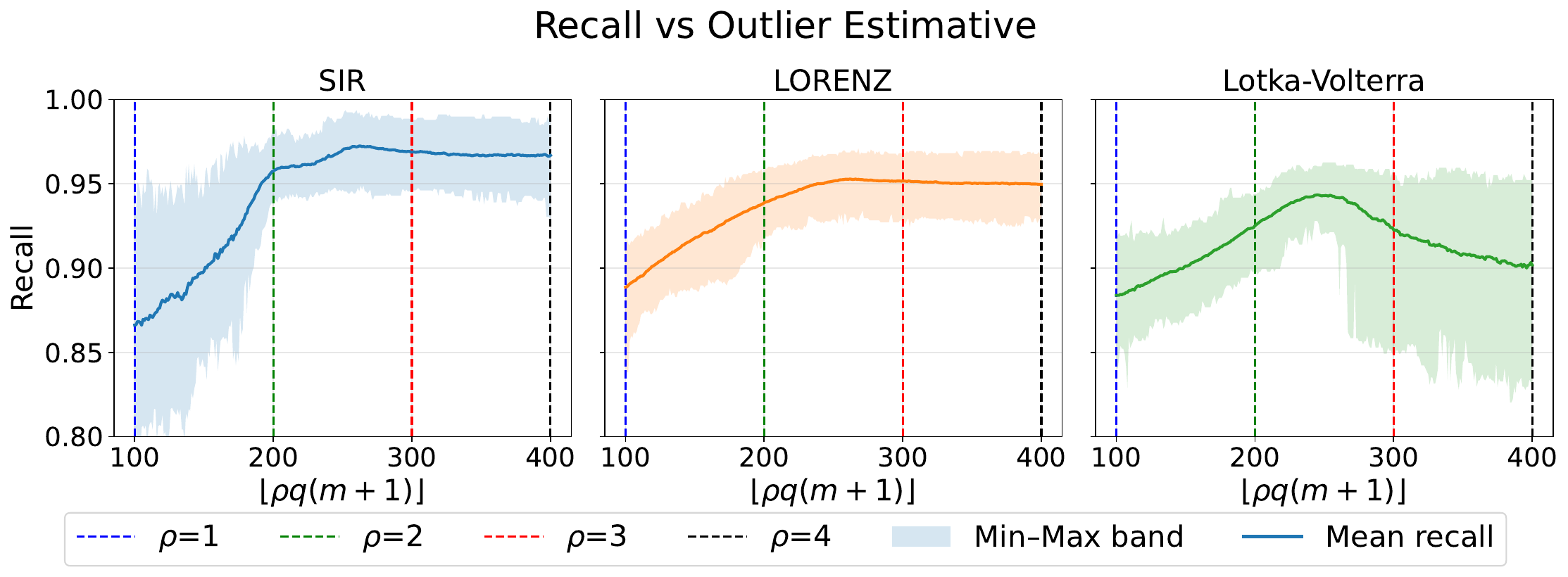}
    \caption{Mean recall of outlier detection using ILTS as a function of the estimated number of outliers. Results are averaged over 100 independent noise realizations with 10\% outliers and a 10\% noise level.}
    \label{fig:q_sweep_10_10}
\end{figure*}

As a result of this configuration, we observe a non-monotonic behavior of the recall as a function of the estimated number of outliers. Figure \ref{fig:q_sweep_10_10} illustrates the mean recall of the outlier detection layer as the outlier estimate used by ILTS is varied. As shown, the estimate that maximizes recall typically lies between $\rho=2$ and $\rho=3$, though the precise optimum may shift depending on the specific noise realization. Leveraging the insights gained from this outlier estimation phase, and considering that an exact system-agnostic estimation is practically unfeasible, we adopt a multiplier factor of $\rho=3$ as a robust, general heuristic rule.

\end{appendices}


\bibliography{references}

\end{document}